\documentclass[12pt,a4paper]{article}

\usepackage[usenames]{color}
\usepackage{hyperref}

\usepackage[centertags]{amsmath}
\usepackage{amsfonts}
\usepackage{amssymb}
\usepackage{amsthm}
\usepackage{latexsym}
\usepackage{newlfont}
\usepackage{verbatim}
\usepackage[english]{babel}
\usepackage{texdraw}
\usepackage{graphicx}
\usepackage[matrix,arrow]{xy}

\newtheorem{thm}{Theorem}[section]
\newtheorem{cor}[thm]{Corollary}
\newtheorem{lem}[thm]{Lemma}
\newtheorem{prop}[thm]{Proposition}

\newtheorem{defn}[thm]{Definition}
\theoremstyle{remark}

\newcommand{\ml}{\textbf}

\newtheorem{lemma}{Lemma}[section]

\theoremstyle{remark}

\newtheorem{example}[lemma]{Example}
\newtheorem{remark}[lemma]{Remark}

\newcommand{\bl}{\begin{lem}}
\newcommand{\el}{\end{lem}}
\newcommand{\bt}{\begin{thm}}
\newcommand{\et}{\end{thm}}
\newcommand{\bc}{\begin{cor}}
\newcommand{\ec}{\end{cor}}
\newcommand{\bp}{\begin{proof}}
\newcommand{\ep}{\end{proof}}
\newcommand{\bpr}{\begin{prop}}
\newcommand{\epr}{\end{prop}}
\newcommand{\brem}{\begin{remark} }
\newcommand{\erem}{\end{remark}}
\newcommand{\bd}{\begin{defn} \em}
\newcommand{\ed}{\end{defn}}
\newcommand{\bex}{\begin{example}
}
\newcommand{\eex}{\end{example}}

\newcommand{\bi}{\begin{itemize}}
\newcommand{\ei}{\end{itemize}}
\newcommand{\ben}{\begin{enumerate} }
\newcommand{\een}{\end{enumerate} }

\newenvironment{enumr}{
\renewcommand{\theenumi}{\roman{enumi}}
\renewcommand{\labelenumi}{$\mathrm{(\theenumi)}$}
\begin{enumerate}     }{\end{enumerate}
\renewcommand{\theenumi}{\arabic{enumi}}
\renewcommand{\labelenumi}{\theenumi.}}

\newcommand{\al}[1]{\forall #1\:}
\newcommand{\exi}[1]{\exists #1\:}

\newlength{\hilflh}

\renewcommand{\emptyset}{\varnothing}

\newcommand{\cL}{{\mathcal L}}

\newcommand{\cP}{{\mathcal P}}

\newcommand{\cM}{{\mathcal M}}

\newcommand{\cC}{{\mathcal C}}

\newcommand{\ga}{\alpha}
\newcommand{\gb}{\beta}
\newcommand{\gd}{\delta}
\renewcommand{\ge}{\varepsilon}
\newcommand{\gl}{\lambda}

\newcommand{\gs}{\sigma}
\newcommand{\gy}{\gamma}
\newcommand{\gw}{\omega}

\renewcommand{\phi}{\varphi}

\newcommand{\eqv}{\leftrightarrow}

\newcommand{\ul}[1]{\underline{#1}}

\renewcommand{\Pr}{\mathsf{Prov}}

\newcommand{\GL}{\mathbf{GL}}
\newcommand{\GLP}{\mathbf{GLP}}

\newcommand{\ch}{\mathrm{ch}}

\newcommand{\Con}{\mathsf{Con}}

\newcommand{\EA}{\mathsf{EA}}

\newcommand{\ZFC}{\mathsf{ZFC}}

\newcommand{\gn}[1]{\ulcorner #1 \urcorner}

\newcommand{\la}{\langle}
\newcommand{\ra}{\rangle}

\renewcommand{\models}{\vDash}      

\newcommand{\nmodels}{\nvDash}

\newcommand{\nCon}{\textsf{-}\mathsf{Con}}
\newcommand{\Imp}{\Rightarrow}

\newcommand{\On}{\mathrm{On}}

\newcommand{\Lim}{\mathrm{Lim}}

\newcommand{\cf}{\mathrm{cf}}

\newcommand{\Reg}{\mathrm{Reg}}

\newcommand{\Log}{\mathrm{Log}}
\newcommand{\tto}{\twoheadrightarrow}

\newcommand{\email}[1]{email: \texttt{#1}}

\begin{document}

\author{Lev Beklemishev\thanks{V.A.~Steklov Mathematical Institute,
RAS; Moscow M.V. Lomonosov State University; National Research
University Higher School of Economics; \email{bekl@mi.ras.ru}}
 \and David
Gabelaia\thanks{TSU Razmadze Institute of Mathematics, Tbilisi; 
\email{gabelaia@gmail.com}}}

\title{Topological interpretations of provability logic}

\maketitle

\abstract{Provability logic concerns the study of modality $\Box$ as
provability in formal systems such as Peano arithmetic. Natural,
albeit quite surprising, topological interpretation of provability
logic has been found in the 1970's by Harold Simmons and Leo Esakia.
They have observed that the dual $\Diamond$ modality, corresponding
to consistency in the context of formal arithmetic, has all the
basic properties of the topological derivative operator acting on a
scattered space. The topic has become a long-term project for the
Georgian school of logic led by Esakia, with occasional
contributions from elsewhere.\newline\indent More recently, a new
impetus came from the study of polymodal provability logic $\GLP$
that was known to be Kripke incomplete and, in general, to have a
more complicated behavior than its unimodal counterpart. Topological
semantics provided a better alternative to Kripke models in the
sense that $\GLP$ was shown to be topologically complete. At the
same time, new fascinating connections with set theory and large
cardinals have emerged.\newline\indent We give a survey of the
results on topological semantics of provability logic starting from
first contributions  by Esakia. However, a special emphasis is put
on the recent work on topological models of polymodal provability
logic. We also included a few results that have not been published
so far, most notably the results of Section \ref{linearity} (due the
second author) and Sections \ref{d-ref}, \ref{ordinal} (due to the
first author).}

\section{Provability algebras and logics}

Provability logics and algebras emerge from, respectively, a modal
logical and an algebraic point of view on the proof-theoretic
phenomena around G\"odel's incompleteness theorems. These theorems
are usually perceived as putting fundamental restrictions on what
can be formally proved in a given axiomatic system (satisfying
modest natural requirements). For the sake of a discussion, we call
a formal theory $T$ \emph{g\"odelian} if \bi
\item $T$ is a first order theory in which the natural numbers along
with the operations $+$ and $\cdot$ are interpretable;
\item $T$ proves some basic properties of these operations and a
modicum of induction (it is sufficient to assume that $T$ contains
the Elementary Arithmetic $\EA$);
\item $T$ has a recursively enumerable set of axioms.
\ei The Second Incompleteness Theorem of Kurt G\"odel (G2) states
that a g\"odelian theory $T$ cannot prove its own consistency
provided it is indeed consistent. More accurately, for any r.e.\
presentation of such a theory $T$, G\"odel has shown how to write
down an arithmetical formula $\Pr_T(x)$ expressing that \emph{$x$ is
(a natural number coding) a formula provable in $T$}. Then the
statement $\Con(T):=\neg \Pr_T(\gn{\bot})$ naturally expresses that
the theory $T$ is consistent. G2 states that $T\nvdash \Con(T)$
provided $T$ is consistent.

Provability logic emerged from the question what properties of
formal provability $\Pr_T$ can be verified in $T$, even if the
consistency of $T$ cannot. Several such
properties have been stated by G\"odel himself \cite{God33}.
Hilbert and Bernays and then L\"ob \cite{Lob55} stated them in the
form of conditions any adequate formalization of a provability
predicate in $T$ must satisfy. After G\"odel's and L\"ob's work it was clear that
the formal provability predicate calls for a treatment as a modality. It led to the
formulation of the G\"odel--L\"ob provability logic $\GL$ and eventually to
the celebrated arithmetical completeness theorem due to Robert Solovay.

Independently, Macintyre and Simmons \cite{MS73} and Magari
\cite{Mag} took a very natural algebraic perspective on the
phenomenon of formal provability which led to the concept of
\emph{diagonalizable algebra}. Such algebras are now more commonly
called \emph{Magari algebras}.
This point of view will be more convenient for our
present purposes.

Recall that the Lindenbaum--Tarski algebra of a theory $T$ is the
set of all $T$-sentences $\text{Sent}_T$ modulo provable equivalence
in $T$, that is, the  structure $\cL_T= \text{Sent}_T/{\sim_T}$
where, for all $\phi,\psi\in\text{Sent}_T$,
\[\phi\sim_T\psi \iff \text{$T\vdash(\phi\eqv\psi)$}.\]
Since we assume $T$ to be based on classical propositional logic,
$\cL_T$ is a boolean algebra with operations $\land$, $\lor$,
$\neg$. Constants $\bot$ and $\top$ are identified with the sets of
refutable and of provable sentences of $T$, respectively. The
standard ordering on $\cL_T$ is defined by
$$[\phi]\leq [\psi] \iff T\vdash \phi\to \psi \iff [\phi\land\psi]=[\phi],$$
where $[\phi]$ denotes the equivalence class of $\phi$.

It is well-known that for consistent g\"odelian $T$ all such
algebras are isomorphic to the unique countable atomless boolean
algebra. (This is a consequence of a strengthening of G\"odel's
First Incompleteness Theorem due to Rosser.) We obtain more
interesting algebras by enriching the structure of the boolean
algebra $\cL_T$ by additional operation(s).

G\"odel's consistency formula induces a unary operator $\Diamond_T$
acting on $\cL_T$:
$$\Diamond_T: [\phi]\longmapsto [\Con(T+\phi)]. $$
The sentence $\Con(T+\phi)$ expressing the consistency of $T$
extended by $\phi$ can be defined as $\neg\Pr_T(\gn{\neg\phi})$. The
dual operator is $\Box_T:[\phi]\longmapsto [\Pr_T(\gn{\phi})],$ thus
$\Box_T x=\neg\Diamond_T\neg x$, for all $x\in\cL_T$.

Bernays--L\"ob derivability conditions ensure that $\Diamond_T$ is
correctly defined on the equivalence classes of the
Lindenbaum--Tarski algebra of $T$. Moreover, it satisfies the
following identities (where we write $\Diamond_T$ simply as
$\Diamond$ and the variables range over arbitrary elements of
$\cL_T$):
\ben
\renewcommand{\labelenumi}{M\theenumi.}
\item $\Diamond \bot = \bot; \quad \Diamond(x\lor y)=\Diamond x\lor \Diamond y$;
\item $\Diamond x=\Diamond (x\land\neg\Diamond x)$.
\een

Notice that Axiom M2 is a formalization of G2 stated for the
theory $T'=T+\phi$, where $[\phi]=x$. In fact, the left hand side
states that $T'$ is consistent, whereas the right hand side states
that $T'+\neg\Con(T')$ is consistent, that is, $T'\nvdash \Con(T')$.
The dual form of Axiom M2, $$\Box(\Box x\to x)=\Box x,$$ expresses
the familiar L\"ob's theorem.

Notice also that M1 induces $\Diamond$ to be monotone: if $x\leq y$ then $\Diamond x \leq \Diamond y$.

Boolean algebras with operators $\cM = (M,\Diamond)$ satisfying the
above identities are called \emph{Magari algebras}. Thus, the main
example of a Magari algebra is the structure $(\cL_T,\Diamond_T)$,
for any consistent g\"odelian theory $T$, which is also called the
\emph{provability algebra of $T$}.

The \emph{transitivity} inequality $\Diamond\Diamond x\leq \Diamond
x$ is often postulated as an additional axiom of Magari algebras,
however it follows from M1 and M2.

\bpr \label{trans} In any Magari algebra $\cM$, there holds
$\Diamond\Diamond x\leq \Diamond x$, for all $x\in M$. \epr

\bp Given any $x\in M$ consider $y:= x\lor \Diamond x$. On the one
hand, we have $$\Diamond \Diamond x \leq (\Diamond
x\lor\Diamond\Diamond x) = \Diamond y.$$ On the other hand, since
$\Diamond x\land\neg\Diamond y=\bot$ we obtain $$\Diamond y\leq
\Diamond (y\land \neg \Diamond y)\leq \Diamond((x\lor\Diamond
x)\land\neg\Diamond y) \leq \Diamond(x\land\neg\Diamond y)\lor
\Diamond \bot\leq \Diamond x.$$ Hence, $\Diamond\Diamond x\leq \Diamond
x$. \ep

In general, we call an \emph{identity} of an algebraic structure
$\cM$ a formula of the form $t(\vec x)=u(\vec x)$, where $t,u$ are
terms, such that $\cM\models\al{\vec x} (t(\vec x)=u(\vec x))$.
Identities of Maragi algebras can be described in terms of modal
logic as follows. Any term (built from the variables using boolean
operations and $\Diamond$) is naturally identified with a formula in
the language of propositional logic with a new unary connective
$\Diamond$. If $\phi(\vec x)$ is such a formula and $\cM$ a Magari
algebra, we write $\cM\models\phi$ iff $\al{\vec x}(t_\phi(\vec
x)=\top)$ is valid in $\cM$, where $t_\phi$ is the term
corresponding to $\phi$. Since any identity in Magari algebras can
be equivalently written in the form $t=\top$, for some term $t$, the
axiomatization of identities of $\cM$ amounts to axiomatizing modal
formulas valid in $\cM$. The \emph{logic of $\cM$}, $\Log(\cM)$, is
the set of all modal formulas valid in $\cM$, that is,
$$\Log(\cM):=\{\phi:\cM\models\phi\},$$ and the logic of a class of
modal algebras is defined similarly.

One of the main parameters of a Magari algebra $\cM$ is its
\emph{characteristic $\ch(\cM)$} which is defined as follows. We let
$$\ch(\cM):= \min \{k\in\gw : \Diamond^k \top = \bot\}$$ and
$\ch(\cM):=\infty,$ if no such $k$ exists. If $T$ is arithmetically
sound, that is, if the arithmetical consequences of $T$ are valid in
the standard model, then $\ch(\cL_T)=\infty$. Theories of finite
characteristics are, in a sense, close to being inconsistent and can
be considered a pathology.

Robert Solovay in his fundamental paper \cite{Sol76} proved that any
identity valid in the structure $(\cL_T,\Diamond_T)$ follows from
the boolean identities together with M1--M2, provided $T$ is
arithmetically sound. This has been generalized by Albert Visser
\cite{Vis84} to arbitrary theories of infinite characteristic. Put
algebraically, this becomes the following statement.

\bt Suppose $\ch(\cL_T,\Diamond_T)=\infty$. An identity holds in
$(\cL_T,\Diamond_T)$ iff it holds in all Magari algebras. \et

Apart from the equational characterization by M1--M2 above, the
identities of Magari algebras can be axiomatized modal-logically. In
fact, the logic of all Magari algebras, and by Solovay theorem the
logic $\Log(\cL_T,\Diamond_T)$ of the Magari algebra of $T$, for any
fixed theory $T$ of infinite characteristic, coincides with the
familiar G\"odel--L\"ob logic $\GL$. By abuse of language we will
often identify $\GL$ with the set of identities of Magari algebras.

A Hilbert-style axiomatization of $\GL$ is usually given in the
modal language where $\Box$ rather than $\Diamond$ is taken as basic
and the latter is treated as an abbreviation for $\neg\Box\neg$. The
axioms and inference rules of $\GL$ are as follows.

\bigskip \textbf{Axiom schemata:} \ben
\renewcommand{\labelenumi}{L\theenumi.}
\item All instances of propositional tautologies;
\item $\Box(\phi\to\psi)\to (\Box\phi\to\Box\psi)$;
\item $\Box(\Box\phi\to \phi)\to \Box\phi$.
\een

\textbf{Rules:} $\phi,\ \phi\to\psi/\psi$ (modus ponens), $\phi/\Box\phi$ (necessitation).

\bigskip
By a well-known result of Segerberg, $\GL$ is sound and complete
w.r.t.\ the class of all transitive and upwards well-founded Kripke
frames. In fact, it is sufficient to restrict the attention to
frames that are finite irreflexive trees. Thus, summarizing various
characterizations above, we have

\bt For any modal formula $\phi$, the following statements are
equivalent:

\begin{enumr} \item $\GL\vdash\phi$;
\item $\phi$ is valid in all Magari algebras;
\item $(\cL_T,\Diamond_T)\models \phi$, for any $T$ of infinite characteristic;
\item $\phi$ is valid in all finite irreflexive tree-like Kripke frames.
\end{enumr}
\et

\section{Topological interpretation}

Having defined the notion of Magari algebra, the first question one
can ask is whether there are any other natural examples of such
algebras apart from the provability algebras of the form
$(\cL_T,\Diamond_T)$. The fact that such algebras naturally emerge
from scattered topological spaces was discovered independently by
Harold Simmons \cite{Sim75} and Leo Esakia
\cite{Esa81}.\footnote{The paper by Simmons was published in the
pre-Solovay era of provability logic and remained largely ignored.
Esakia apparently independently discovered the interpretation when
he was already familiar with the works of Magari and Solovay.} We
now infer this semantics from rather general considerations.

It is a standard idea in logic, going back at least as early as to
Leibnitz, to interpret propositions as subsets of a given nonempty
set $X$. Then, logical connectives are naturally translated as the
familiar boolean operations on sets. From this point of view,
$\Diamond$ should correspond to an operator acting on the subsets of
$X$, so we come to the following mathematical question.

Let $X$ be a nonempty set, $\cP(X)$ the boolean algebra of subsets
of $X$. Consider any operator $\delta:\cP(X)\to \cP(X)$ and the
structure $(\cP(X),\delta)$. Can $(\cP(X),\delta)$ be a Magari
algebra and, if yes, when? Simmons and Esakia, in some sense, found
a canonical answer to this question.

Let $(X,\tau)$ be a topological space, where $\tau$ denotes the set
of open subsets of $X$, and let $A\subseteq X$. Topological
\emph{derivative $d_\tau(A)$ of $A$} is the set of limit points of
$A$:
\[
x\in d_\tau (A) \iff \forall U\in\tau\:(x\in U\Imp \exists y\neq x\:
(y\in U\cap A)).
\]
Notice that $c_\tau (A):=A\cup d_\tau(A)$ is the
closure of $A$ and $i_\tau(A):=A \setminus d_\tau(A)$ is the set of isolated points of
$A$.

The classical notion of a scattered topological space is due to
Georg Cantor. $(X,\tau)$ is called \emph{scattered} if every
nonempty subspace $A\subseteq X$ has an isolated point.

\bt[Simmons, Esakia] \label{sim-esa} The following statements are
equivalent:
\begin{enumr}
\item $(X,\tau)$ is scattered;
\item For all $A\subseteq X$, $d_\tau(A)=d_\tau(A\setminus d_\tau(A))$;
\item $(\cP(X),d_\tau)$ is a Magari algebra.
\end{enumr}
\et Notice that (ii) means that each point of $A$ is a limit point
of its isolated points. Algebras of the form $(\cP(X),d_\tau)$
associated with a topological space $(X,\tau)$ will be called
\emph{derivative algebras}. Thus, this theorem states that the
derivative algebra of $(X,\tau)$ is Magari iff $(X,\tau)$ is
scattered.

\bp First, we notice that Statement (iii) is just a reformulation of
(ii). In fact, the identities $d_\tau\emptyset=\emptyset$ and
$d_\tau(A\cup B)=d_\tau (A)\cup d_\tau (B)$ are obviously valid in
all topological spaces. We show the equivalence of (i) and (ii).

Suppose $(X,\tau)$ is scattered, $A\subseteq X$ and $x\in
d_\tau(A)$. Consider any open neighborhood $U$ of $x$. Since $U\cap
A\setminus \{x\}$ is nonempty, it has an isolated point $y\neq x$. Since $U$
is open, $y$ is an isolated point of $A$, that is, $y\in
A\setminus d_\tau(A)$. Hence, $x\in d_\tau(A\setminus d_\tau(A))$. The inclusion
$d_\tau(A\setminus d_\tau(A))\subseteq d_\tau(A)$ follows from the
monotonicity of $d_\tau$. Therefore Statement (ii) holds.

Suppose that (ii) holds and let $A\subseteq X$ be nonempty. We show
that $A$ has an isolated point. If $d_\tau A$ is empty, we are done.
Otherwise, take any $x\in d_\tau A$. Since $x$ is a limit of
isolated points of $A$, there must be at least one such point. \ep

We notice that the transitivity principle $d_\tau d_\tau A\subseteq
d_\tau A$ topologically means that the set $d_\tau A$, for any
$A\subseteq X$, is closed. We recall the following standard equivalent
characterization an easy proof of which we shall omit.

\bpr For any topological space $(X,\tau)$, the following statements
are equivalent:
\begin{enumr}
\item Every $x\in X$ is an intersection of an open and a closed set;
\item For each $A\subseteq X$, the set $d_\tau A$ is closed.
\end{enumr}
\epr

Topological spaces satisfying either of these conditions are called
\emph{$T_d$-spaces}. Condition (i) shows that $T_d$ is a weak
separation property located between $T_0$ and $T_1$. Thus,
Proposition \ref{trans} yields the following corollary.

\bc All scattered spaces are $T_d$. \ec

We have seen in Theorem \ref{sim-esa} that  each scattered space
equipped with a topological derivative operator is a Magari algebra.
The following result shows that any Magari algebra on $\cP(X)$ can
be described in this way.

\bt \label{mag-scat} If $(\cP(X),\delta)$ is a Magari algebra, then
$X$ bears a unique topology $\tau$ for which
$\delta=d_\tau$. Moreover, $\tau$ is scattered. \et

\bp

We first remark that, if $(\cP(X),\gd)$ is a Magari algebra, then
the operator $c(A):= A\cup \gd A$ satisfies the Kuratowski axioms of
the topological closure: $c\emptyset=\emptyset$, $c(A\cup B)=cA\cup
cB$, $A\subseteq cA$, $ccA=cA$ (the latter uses the transitivity
property of $\gd$). This defines a topology $\tau$ on $X$ in which a
set $A$ is closed iff $A=c(A)$ iff $\gd A\subseteq A$. If $\gd$ coincides
with $d_\tau$, this condition is also clearly necessary. So, the required topology
is unique. To show that $\gd=d_\tau$ indeed, we
need an auxiliary lemma.

\bl Suppose $(\cP(X),\gd)$ is Magari. Then, for all $x\in X$,
\begin{enumr}
\item $x\notin \gd(\{x\})$;
\item $x\in \gd A \iff x\in \gd(A\setminus\{x\})$.
\end{enumr}
\el

\bp (i) By Axiom M2 we have $\gd\{x\}\subseteq \gd(\{x\}\setminus
\gd\{x\})$. If $x\in \gd\{x\}$ then $\gd(\{x\}\setminus
\gd\{x\})=\gd\emptyset=\emptyset$.
Hence, $\gd\{x\}=\emptyset$, a contradiction.

(ii) $x\in \gd A$ implies $x\in \gd ((A\setminus\{x\})\cup \{x\})=
\gd(A\setminus\{x\})\cup \gd\{x\}$. By (i), $x\notin\gd\{x\}$, hence
$x\in \gd(A\setminus\{x\})$. The other implication follows from the
monotonicity of $\gd$. \ep

\bl Suppose $(\cP(X),\gd)$ is Magari and $\tau$ is the associated
topology. Then $\gd=d_\tau$. \el

\bp Let $d=d_\tau$; we show that, for any set $A\subseteq X$,
$dA=\gd A$. Notice that, for any $B$, $cB=dB\cup B=\gd B\cup B$.

Assume $x\in \gd A$ then $$x\in \gd(A\setminus\{x\})\subseteq
c(A\setminus\{x\})\subseteq d(A\setminus \{x\})\cup
(A\setminus\{x\}).$$ Since $x\notin A\setminus\{x\}$, we obtain
$x\in d(A\setminus \{x\})$. By the monotonicity of $d$, $x\in dA $.
Similarly, if $x\in dA$ then $x\in d(A\setminus\{x\})$. Hence,
$$x\in c(A\setminus\{x\})=\gd(A\setminus\{x\})\cup (A\setminus
\{x\}).$$ Since $x\notin A\setminus\{x\}$ we obtain $x\in \gd A$.
\ep

From this lemma and Theorem \ref{sim-esa} we infer that the
associated topology $\tau$ must be scattered. \ep

Theorem \ref{mag-scat} shows that to study a natural set-theoretic
interpretation of provability logic means to study the semantics
of $\Diamond$ as a derivative operation on a scattered topological space.

Derivative semantics of modality was first suggested in the fundamental paper by McKinsey and Tarski \cite{MT44}. See the paper by Kudinov and Shehtman in this volume for a detailed survey of such semantics for general topological spaces. The emphasis in this paper is on the logics related to formal provability. Specifically, scattered spaces (and their generalizations) will play a major role here.

\section{Topological completeness theorems} \label{topcompl}

Natural examples of scattered topological spaces come from
orderings. Two examples will play an important role below.

Let $(X,\prec)$ be a strict partial ordering. The \emph{left
topology $\tau_{\leftarrow}$ on $(X,\prec)$} is given by all sets
$A\subseteq X$ such that $\forall x,y \:( y\prec x\in A \Imp y\in
A)$. We obviously have that $(X,\prec)$ is well-founded iff
$(X,\tau_{\leftarrow})$ is scattered. The \emph{right topology} or
the \emph{upset topology} is defined similarly.

This topology is, in general, non-Hausdorff. More natural is the
\emph{order topology} on a linear ordering $(X,<)$, which is
generated by all intervals $(\ga,\gb)$ such that $\ga,\gb\in
X\cup\{\pm\infty\}$ and $\ga<\gb$. The order topology refines the
left topology and is scattered on any ordinal.

Given a topological space $(X,\tau)$, we shall
denote the logic of its derivative algebra $(\cP(X),d_\tau)$ by
$\Log(X,\tau)$, and we use similar notation for the logic of a class
of derivative algebras.

The logic $\Log(X,\tau)$ is often equivalently defined in terms of
valuations as follows. A \emph{valuation on $(X,\tau)$} is any map
$v$ from the set of propositional formulas to the powerset of $X$
satisfying the following conditions: \bi
\item $v(\phi\lor \psi)=v(\phi)\cup v(\psi)$, $v(\neg \phi)= X\setminus v(\phi)$, $v(\top)=X$, $v(\bot)=\emptyset$,
\item $v(\Diamond\phi) = d(v(\phi))$.
\ei One usually writes $X,x\models_v\phi$ for $x\in v(\phi)$. A
formula $\phi$ is \emph{valid in $X$} (denoted $X\models\phi$) if
$\forall v, v(\phi)=X$. Then, $\Log(X,\tau)$ coincides with the set
of all formulas valid in $(X,\tau)$. Similarly, if $\cC$ is a class
of spaces, $\Log(\cC)$ is the set of formulas valid in all members
$(X,\tau)\in\cC$.

Since $\Log(\cC)$ is the logic of a class of Magari algebras, this
set of formulas is closed under the rules modus ponens,
necessitation and substitution, and contains all theorems of
G\"odel--L\"ob logic $\GL$. Thus, $\Log(\cC)$ is a normal modal
logic extending $\GL$, for any class $\cC$ of scattered spaces.

Topological interpretation of provability logic suggests a question
whether $\GL$ is complete for any natural class of scattered
topological spaces. Leo Esakia \cite{Esa81} has shown that $\GL$
coincides with the logic of the class of all scattered topological spaces. In fact, it is sufficient to take finite trees (considered as strict partial orderings) with the associated upset topology.

\bt[Esakia]
\begin{enumr}
\item $\Log(\cC)=\GL$, where $\cC$ is the class of all finite irreflexive trees taken with the upset topology.
\item There is a countable scattered space $(X,\tau)$ such that $\Log(X,\tau)=\GL$.
\end{enumr}
\et

We notice that this theorem follows almost immediately from the
completeness of $\GL$ w.r.t.\ its Kripke semantics. In fact, if
$(X,\prec)$ is a strict partial ordering, then the modal algebra
associated with the Kripke frame $(X,\prec)$ is the same as the
derivative algebra of $(X,\tau)$ where $\tau$ is its upset topology.
This implies that any modal logic of a class of strict partial
orders is complete w.r.t.\ topological semantics.

Esakia's theorem has been strengthened by Abashidze~\cite{Aba85} and
Blass~\cite{Bla90}.

\bt[Abashidze, Blass] Consider an ordinal $\Omega\geq \gw^\gw$
equipped with the order topology. Then $\Log(\Omega)=\GL$. \et

Thus, $\GL$ is complete w.r.t.\ a natural scattered topological
space. In Section \ref{s:ab-blass} we give an alternative proof of
this result by utilizing the recursive constructions of finite
irreflexive trees and ordinals below $\gw^\gw$.

\section{Some useful techniques}

\subsection{Ranks and $d$-maps}

An equivalent characterization of scattered spaces is often given in
terms of the following transfinite \emph{Cantor--Bendixson sequence}
of subsets of a topological space $(X,\tau)$:
\begin{itemize}
\item $d_\tau^0 X=X$;
\quad $d_\tau^{\alpha+1} X=d_\tau(d_\tau^{\alpha} X)$ and
\item $d_\tau^\alpha X=\bigcap\limits_{\beta<\alpha}d_\tau^\beta X$ if $\alpha$ is a limit
ordinal.
\end{itemize}


It is easy to show by transfinite induction that, for any
$(X,\tau)$, all sets $d_\tau^\ga X$ are closed and that $d_\tau^\ga
X \supseteq d_\tau^\gb X$ whenever $\ga\leq \gb$.

\bt[Cantor] $(X,\tau)$ is scattered iff $d_\tau^\ga X = \emptyset$,
for some ordinal $\ga$. \et

\bp Let $d=d_\tau$. If $(X,\tau)$ is scattered then we have $d^\ga
X\supset d^{\ga+1}X$, for each $\ga$ such that $d^\ga
X\neq\emptyset$. By cardinality considerations this yields an $\ga$
such that $d^\ga X=\emptyset$.

Conversely, suppose $A\subseteq X$ is nonempty. Let $\ga$ be the
least ordinal such that $A\nsubseteq d^\ga X$. Obviously, $\ga$
cannot be a limit ordinal, hence $\ga=\gb+1$ for some $\gb$ and
there is an $x\in A \setminus d^{\gb+1}X$. Since $A\subseteq d^\gb
X$, we also have $x\in d^\gb X$. Since $x\notin d^{\gb+1}X=d(d^\gb
X)$, $x$ is isolated in the relative topology of $d^\gb X$, and
hence in the relative topology of $A\subseteq d^\gb X$. \ep

Call the least $\ga$ such that $d_\tau^\ga X=\emptyset$ the
\emph{Cantor--Bendixson rank} of $X$ and denote it by
$\rho_\tau(X)$. The \emph{rank function} $\rho_\tau:X\to \On$ is
defined by
$$\rho_\tau(x):=\min\{\ga: x\notin d_\tau^{\ga+1}(X)\}.$$
Notice that $\rho_\tau$ maps $X$ onto
$\rho_\tau(X)=\{\ga:\ga<\rho_\tau(X)\}$. Also, $\rho_\tau(x)\geq\ga$
iff $x\in d_\tau^\ga X$. We omit the subscript $\tau$ whenever there
is no danger of confusion.

\bex Let $\Omega$ be an ordinal equipped with its \emph{left
topology}. Then $\rho(\ga)=\ga$, for all $\ga$. \eex

\bex Let $\Omega$ be an ordinal equipped with its \emph{order
topology}. Then $\rho$ is the function $\ell$ defined by
$$\ell(0)=0; \quad \ell(\ga)=\gb \text{ if $\ga=\gy+\gw^{\gb}$, for some
$\gy$, $\gb$.}$$ By the Cantor normal form theorem, for any $\ga>0$,
such a $\gb$ is uniquely defined. \eex

Let $(X,\tau_X)$ and $(Y,\tau_Y)$ be topological spaces, and let
$d_X$, $d_Y$ denote the corresponding derivative operators. A map
$f:X\to Y$ is called a \emph{d-map} if $f$ is continuous, open and
\emph{pointwise discrete}, that is, $f^{-1}(y)$ is a discrete
subspace of $X$ for each $y\in Y$. $d$-maps are well-known to
satisfy the properties expressed in the following lemma (see
\cite{BEG05}).

\bl\label{l:d-maps}
\begin{enumr}
\item $f^{-1}(d_{Y}(A))= d_{X}(f^{-1}(A))$, for any $A\subseteq Y$;
\item $f^{-1}:(\cP(Y),d_{Y})\to (\cP(X),d_{X})$ is a homomorphism of derivative algebras;
\item If $f$ is onto, then $\Log(X,\tau_X)\subseteq \Log(Y,\tau_Y)$.
\end{enumr}
\el

In fact, (i) is easy to check directly; (ii) follows from (i) and
(iii) follows from (ii). Each of the conditions (i) and
(ii) is equivalent to $f$ being a $d$-map.

From lemma~\ref{l:d-maps}(i) we easily obtain the following
corollary by transfinite induction.

\bc\label{CB-pres} Suppose $f:X\to Y$ is a $d$-map. Then, for each
ordinal $\ga$, $d_{X}^\ga X= f^{-1}(d_{Y}^\ga Y)$. \ec

The following lemma states that the rank function, when the ordinals
are equipped with their left topology, becomes a $d$-map. It is also
uniquely characterized by this property.

\bl \label{rank-dmap} Let $\Omega$ be the ordinal $\rho_\tau(X)$
taken with its left topology. Then
\begin{enumr}
\item $\rho_\tau:X\tto \Omega$ is an onto $d$-map;
\item If $f:X\to \gl$ is a $d$-map, where $\gl$ is an ordinal with its left topology, then $f(X)=\Omega$ and $f=\rho_\tau$.
\end{enumr}
\el

\bp Let $\rho$ denote $\rho_\tau$.

(i) $\rho$ is continuous, because the set
$\rho^{-1}[0,\ga)=X\setminus d^{\ga}X$ is open.

$\rho$ being open means that, for each open $U\subseteq X$, whenever
$\ga \in \rho(U)$ and $\gb<\ga$ one has $\gb\in \rho(U)$. Fix an
$x\in U$ such that $\rho(x)=\ga$. Consider the set
$X_\gb:=\rho^{-1}(\gb)=d^\gb X\setminus d(d^\gb X)$. For any subset
$A$ of a scattered space we have $d(A)=d(A\setminus dA)$, hence $d
X_\gb=d(d^\gb X)\subseteq d^\ga X$. Since $\rho(x)=\ga$ it follows
that $x\in d X_\gb$. Hence $U\cap X_\gb\neq\emptyset$, that is,
$\gb\in \rho(U)$.

$\rho$ being pointwise discrete means $X_\ga=\rho^{-1}(\ga)$ is
discrete, for each $\ga$. In fact, $X_\ga=d^\ga X\setminus d(d^{\ga}
X)$ is the set of isolated points of $d^\ga X$. Thus, it cannot help
being  discrete.

(ii) Since $f$ is a $d$-map, by Corollary \ref{CB-pres} we obtain
that $f^{-1}[\ga,\gl)=d^\ga X$, for each $\ga<\gl$. Hence,
$f^{-1}(\ga)=\rho^{-1}(\ga)$, for each $\ga<\gl$, that is, $f=\rho$
and $f(X)=\rho(X)=\Omega$. \ep

Let $\rho_X$ and $\rho_Y$ denote the rank functions of $(X,\tau_X)$
and $(Y,\tau_Y)$, respectively.

\bc\label{c:ranks compose}  If $f:X\to Y$ is a $d$-map, then $\rho_{X}=\rho_{Y}\circ f$.
\ec

\bp Clearly, $\rho_{Y}\circ f: X\to \Omega$ is a $d$-map. Statement
(ii) of the previous lemma yields the result. \ep

\subsection{The d-sum construction}

Given a tree $T$, one can build a new tree by `plugging in' other
trees in place of the maximal points (the `leaves') of $T$.
Similarly, given an ordinal $\alpha$, one can `plug in' new ordinals
$\alpha_i$ for each isolated (i.e. non-limit) point $i\in\alpha$ to
obtain another ordinal. These constructions of `summing up' spaces
or relational structures `along' another space or a relational
structure are a special case of a general construction, which proved
rather useful in its various manifestations for proving topological
completeness theorems. Here we present a variant of this
construction, called \emph{$d$-sum}, which can be used to
recursively build both finite irreflexive trees and ordinals. We will point out a particular application of
$d$-sums by sketching a proof of Abashidze--Blass theorem. For a
more substantial application of the $d$-sum construction, in which
the summed spaces are homeomorphic to each other, see
\cite{BekGab11} (where the corresponding construction is called
$d$-product).


\bd Let $X$ be a topological space and let $\{Y_j\mid j\in i(X)\}$
be a collection of spaces indexed by the set $i(X)$ of isolated
points of $X$. We uniquely extend it to the collection $\{Y_j\mid
j\in X\}$ by letting $Y_j=\{j\}$ for all $j\in d_\tau X$.

We define the \emph{d-sum} $(Z,\tau_Z)$ of $\{Y_j\}$ over $X$
(denoted $\sum^d_{j\in X} Y_j$) as follows. The base set is the
disjoint union $Z:=\bigsqcup_{j\in X} Y_j$. Define the map $\pi:Z\to
X$ by putting $\pi(y)=j$ whenever $y\in Y_j$. Now let the topology
$\tau_Z$ consist of the sets $V\cup\pi^{-1}(U)$ where $V$ is open in
the topological sum $\bigsqcup_{j\in i(X)} Y_j$ and $U$ is open in
$X$. It is not difficult to check that $\tau_Z$ qualifies for a
topology. \ed


As an immediate application of this construction we present a way of
building finite irreflexive trees (called simply \emph{trees} for
brevity henceforth) recursively using the
$d$-sum construction.

\bd ($n$-fork) \\
Let $\mathfrak F_n=(W_n,R_n)$ be defined as follows:
$W_n=\{r,w_0,w_1,\dots,w_{n-1}\}$ and $R_n=\{(r,w_i)\mid 0\leq
i<n\}$. We will consider $\mathfrak F_n$ equipped with the upset
topology w.r.t.\ $R_n$. \ed

We view trees, and $n$-forks in particular, as topological spaces.

\begin{defn}(trees)
\begin{enumerate}
\item[(t1)] An irreflexive point is a tree.
\item[(t2)] An $n$-fork $\mathfrak F_n$ is a tree.
\item[(t3)] The $d$-sum of trees is a tree.
\end{enumerate}
\end{defn}

It is straightforward (if somewhat tedious) to check that trees thus
obtained are precisely the finite irreflexive trees. Note that the leaves of a tree are the isolated points in the topology. Therefore, taking the $d$-sum of trees $T_i$ over a tree $T$ means simply plugging in $T_i$'s in place of the leaves of $T$. It is worth noting that a tree might be obtained in various different ways as a d-sum of other trees.

In a similar fashion the $d$-sum construction can be applied to ordinals to obtain larger ordinals. For example, it is easy to see that summing up ordinals $\omega^n+1$ along the ordinal $\omega +1$ produces the ordinal $\omega^{n+1}+1$. In general, from the well-known theorem characterizing countable ordinals as countable locally compact scattered Hausdorff spaces, it follows that the $d$-sum of countable ordinals is a countable ordinal. In fact, even more general statement is true. Namely, if one takes a family of ordinals $(\alpha_i)_{i\in \beta}$ such that $\alpha_i = 1$ for limit $i$, then the ordinal sum $\sum_{i\in\beta} \alpha_i$ is an ordinal that is homeomorphic to the d-sum $\sum^d_{i\in\beta} \alpha_i$. This can be checked directly, by examining the descriptions of neighborhoods in respective topologies. It follows that the $d$-sum of ordinals is an ordinal.

The last paragraph reveals the intrinsic similarity between the inner structure of trees and ordinals. To fix this similarity, we employ the following general lemma showing that $d$-sums, in a way, commute with $d$-maps. More precisely:

\bl\label{t:d-sum d-maps to d-sum} Let $X$ and $X'$ be two spaces
and let $\{Y_j\mid j\in i(X)\}$ and $\{Y'_k\mid k\in i(X')\}$ be
collections of spaces indexed by $i(X)$ and $i(X')$, respectively.
Suppose further that $f:X\to X'$ is an onto $d$-map, and for each
$j\in i(X)$ there is an onto $d$-map $f_j:Y_j\to Y'_{f(j)}$. Then
there exists an onto $d$-map $g:\sum^d_{j\in X} Y_j\to \sum^d_{k\in
X'} Y'_k$. \el

\bp First note that since $f$ is a $d$-map, $f(j)$ is isolated in
$X'$ iff $j$ is isolated in $X$. Indeed, by openness of $f$, if
$\{j\}\in\tau$, then $\{f(j)\}\in\tau'$. Conversely, if $f(j)$ is
isolated, then $f^{-1}f(j)$ is both open and discrete, by continuity
and pointwise discreteness of $f$. Hence, any point in $f^{-1}f(j)$,
and $j$ in particular, is isolated in $X$. For convenience, let us
denote $f_*\equiv f{\upharpoonright}_{d_\tau X}$ and $f^*\equiv
f{\upharpoonright}_{iX}$. It follows that $f^*: i(X)\to i(X')$ and
$f_*: d_\tau X\to d_{\tau'} X'$ are well-defined onto maps and
$f=f^*\cup f_*$. Thus, in particular, the space $Y'_{f(j)}$ in the
formulation of the theorem is well-defined.

Take $g$ to be the set-theoretic union $g=f_*\cup\bigcup_{j\in iX}
f_j$. We show that $g$ is a $d$-map. Let $\pi$ and $\pi'$ be the
`projection' maps associated with $\sum^d_{j\in X} Y_j$ and
$\sum^d_{k\in X'} Y'_k$, respectively. To show that $g$ is open,
take $W=V\cup \pi^{-1}(U)\in\tau_Z$. Then $g(W)=g(V)\cup
g(\pi^{-1}(U))$. That $g(V)$ is open in the topological sum of
$Y'_k$ is clear from the openness of the maps $f_j$. Moreover, from
the definition of $g$ and the fact that all $f_j$ are onto it can be
easily deduced that $g(\pi^{-1}(U))=\pi'^{-1}(f(U))$. Since $f$ is
an open map, it follows that $g(W)$ is open in $\tau'_Z$. To see
that $g$ is continuous, take $W'=V'\cup \pi'^{-1}(U')\in\tau'_Z$.
Then $g^{-1}(W')=g^{-1}(U')\cup g^{-1}(\pi'^{-1}(U'))$. Again, the
openness of $g^{-1}(U')$ is trivial. It is also easily seen that
$g^{-1}(\pi'^{-1}(U'))=\pi^{-1}(f^{-1}(U'))$. It follows that
$g^{-1}(W')$ is open in $\tau_Z$. To see that $g$ is pointwise
discrete is straightforward, given that $f$ and all the $f_j$ are
pointwise discrete. \ep

In the next section we will use this lemma to prove the
Abashidze--Blass theorem.

\section{Proof of Abashidze--Blass theorem}\label{s:ab-blass}


The topological completeness theorem proved independently by
Abashidze and Blass establishes that each non-theorem of $\ml{GL}$
can be refuted on any ordinal greater than $\omega^\omega$. A
simplified proof of this result appeared recently in \cite[Theorem
3.5]{BM10}. The crucial part of that proof rests on showing that
each tree of depth $n$ is an onto image of an ordinal $\omega^n+1$
via a $d$-map \cite[Lemma 3.4]{BM10}. That is probably the simplest
and most explicit proof of Abashidze--Blass theorem, however by our
approach we aim to illuminate the underlying recursive mechanism,
which might offer additional insight. Here is an analog of
\cite[Lemma 3.4]{BM10} proved using the $d$-sum construction:

\bl \label{t:trees are d-maps of s-ordinals} For each tree $T$ there
exists a countable ordinal $\alpha<\omega^\omega$ and an onto map $f:\alpha\to T$ such that
$f$ is a $d$-map. \el
\bp [Sketch] The proof proceeds by structural
induction. It is clear that the claim is true for a one-point tree.
Let us consider now an $n$-fork $\mathfrak F_n$ and build a $d$-map
$f$ from $\omega+1$ onto $\mathfrak F_n$. Let $f(x)=w_{x \bmod n}$
for $x<\omega$ and let $f(\omega)=r$. It is straightforward to check
that $f$ is a $d$-map. That each tree is an image of an appropriate
ordinal via a $d$-map now follows from lemma~\ref{t:d-sum d-maps to d-sum} using induction on the depth of trees. Indeed, suppose the theorem is true for all trees of depth less than $n$. Take a tree $T$ of depth $n$. It is clear that $T$ can be presented as a $d$-sum of trees of strictly smaller depth in various ways. Using the induction hypothesis, each of the smaller trees is an image of a countable ordinal under a $d$-map. Applying lemma~\ref{t:d-sum d-maps to d-sum} and observing that the $d$-sum of countable ordinals is a countable ordinal produces a countable ordinal $\alpha$ and an onto d-map $f:\alpha\to T$. By corollary~\ref{c:ranks compose} the rank of $\alpha$ is equal to the rank of $T$, that is, to $n$. It follows that $\alpha<\omega^\omega$, which finishes the proof.\ep

One of the ways to present the tree $T$ as a $d$-sum of smaller trees is to `prune' $T$ immediately after the root. Then $T$ is seen as a $d$-sum of smaller trees along a $k$-fork for some finite $k$. This is precisely the approach usually taken in proving the above lemma. However, the $d$-sum construction allows for greater flexibility. For example, $T$ can be pruned immediately below the leaves, obtaining finitely many forks and a tree of depth $n-1$. Or indeed, $T$ can be pruned \emph{anywhere}, provided each of the branches of maximal length $n$ is pruned.

And finally, the Abashidze--Blass theorem:

\bt[Abashidze--Blass] The logic of any ordinal
$\alpha\geq\omega^\omega$ is $\ml{GL}$. \et \bp Take a non-theorem
$\varphi$ of $\ml{GL}$. Then $\varphi$ can be refuted on a finite
tree $T$. By Lemma~\ref{t:trees are d-maps of s-ordinals} there
exists an ordinal $\beta<\omega^\omega$ that maps onto $T$ via a
$d$-map. Since $d$-maps reflect satisfiability (the contrapositive
of Lemma~\ref{l:d-maps}(iii)), $\varphi$ can be refuted on $\beta$.
But $\beta$ is an open subspace of $\alpha$. It follows that
$\varphi$ can be refuted on $\alpha$. \ep


\section{Topological semantics of linearity axioms} \label{linearity}

For a g\"odelian theory $T$ consider the subalgebra $\cL_T^0$ of
$(\cL_T,\Diamond_T)$ generated by $\top$. If
$\ch(\cL_T,\Diamond_T)=\infty$, then, clearly,
$\ch(\cL_T^0,\Diamond_T)=\infty$. In fact, the Magari algebra
$(\cL_T^0,\Diamond_T)$ is known (see \cite{IcaJoo12}) to have the
logic $\ml{GL.3}$ which is obtained from $\ml{GL}$ by adding the
following axiom:
 \[
(.3)\ \ \ \  \Diamond p\wedge \Diamond q\to \Diamond(p\wedge
q)\vee\Diamond(p\wedge \Diamond q)\vee\Diamond(\Diamond p\wedge q)
 \]
This is a so called `linearity axiom' and, as the name suggests, its
finite rooted Kripke frames are precisely the finite strict linear
orders. Since $\ml{GL.3}$ is Kripke complete \cite{ChZa}, its
topological completeness is almost immediate. However, it is not
immediately clear what kind of scattered spaces does the linearity
axiom isolate. To characterize GL.3-spaces, let us first simplify
the axiom $(.3)$. Let $(lin)$ be the following formula:
\[
\Box(\Box^+p\vee\Box^+q)\to\Box p\vee\Box q
\]
where $\Box^+\varphi$ is a shorthand for $\varphi\wedge\Box\varphi$.

\bl In $\ml{GL}$ the schema $(.3)$ is equivalent to $(lin)$.\el \bp
To show that $(lin)\vdash_{\ml{GL}}(.3)$, witness the following
syntactic argument. Observe that the dual form of $(lin)$ looks as
follows:
\[
\Diamond p\wedge \Diamond
q\to\Diamond(\Diamond^+p\wedge\Diamond^+q)\ \ \ \ \ \ \ \ \ \ \ (*)
\]
where $\Diamond^+\varphi:=\varphi\vee\Diamond\varphi$.
Furthermore, an instance of the $\ml{GL}$ axiom looks as follows:
\[
\Diamond(\Diamond^+p\wedge\Diamond^+q)\to\Diamond(\Diamond^+p\wedge\Diamond^+q\wedge\Box(\Box^+\neg
p\vee \Box^+\neg q)).
\]
By the axiom $(lin)$ we also have:
\[
\Box(\Box^+\neg p\vee \Box^+\neg q)\to(\Box\neg p\vee\Box\neg q).
\]
So, using the monotonicity of $\Diamond$ we obtain:
\[
\Diamond p\wedge \Diamond
q\to\Diamond(\Diamond^+p\wedge\Diamond^+q\wedge(\Box\neg p\vee
\Box\neg q)).
\]
By Boolean logic we clearly have:
\[\Diamond^+p\wedge\Diamond^+q \eqv (p\wedge q)\vee(p\wedge \Diamond q)\vee(\Diamond p\wedge q)\vee (\Diamond p\wedge \Diamond q)\ \ \ \ \ \ (**)
\]
and
\[
(\Box\neg p\vee \Box\neg q)\eqv \neg(\Diamond p\wedge\Diamond q).
\]
Using these, together with the monotonicity of $\Diamond$ we finally
arrive at:
\[
\Diamond p\wedge \Diamond q\to\Diamond((p\wedge q)\vee(p\wedge
\Diamond q)\vee(\Diamond p\wedge q)),
\]
which is equivalent to $(.3)$ since $\Diamond$ distributes over
$\vee$.

To show the converse, we observe that $(.3)$ implies $(lin)$ even in
the system $\ml{K}$. Indeed, the formula $(*)$, which is the dual
form of $(lin)$, can be rewritten, using $(**)$ and the distribution
of $\Diamond$ over $\vee$ as follows:
\[
\Diamond p\wedge \Diamond q\to\Diamond(p\wedge
q)\vee\Diamond(p\wedge \Diamond q)\vee\Diamond(\Diamond p\wedge
q)\vee \Diamond(\Diamond p\wedge \Diamond q),
\]
which is clearly a weakening of $(.3)$. Therefore
$(.3)\vdash_{\ml{GL}}(lin)$.\ep

It follows that a scattered space is a GL.3-space iff it validates
$(lin)$. To characterize such spaces, consider the following
definition.

\bd Call a scattered space \emph{primal} if for each $x\in X$ and
$U,V\in\tau$, $\{x\}\cup U\cup V\in\tau$ implies $\{x\}\cup
U\in\tau$ or $\{x\}\cup V\in\tau$. \ed

In other words, $X$ is primal if the collection of punctured
neighborhoods of each non-isolated point is a prime filter in the
Heyting algebra $\tau$.
\bt Let $X$ be a scattered space. Then
$X\models(lin)$ iff $X$ is primal.\et
\bp Let $X$ be a scattered
space together with a valuation $v$. Let $P:=v(p)$ and $Q:=v(q)$ denote the
truth-sets of $p$ and $q$ respectively. Then the truth sets of
$\Box^+p$ and $\Box^+q$ are $I_\tau P$ and $I_\tau Q$, where
$I_\tau$ is the interior operator of $X$. We write $x\models\phi$ for $X,x\models_v\phi$ (as introduced in section~\ref{topcompl}).

Suppose $X$ is primal and for some valuation
$x\models\Box(\Box^+p\vee\Box^+q)$. Then there exists an open
neighborhood $W$ of $x$ such that
$W\setminus\{x\}\models\Box^+p\vee\Box^+q$. In other words,
$W\setminus\{x\}\subseteq I_\tau P\cup I_\tau Q$. Let $U=W\cap
I_\tau P\in\tau$ and $V=W\cap I_\tau Q\in\tau$. Then $\{x\}\cup
U\cup V=W\in\tau$. It follows that either $\{x\}\cup U\in\tau$ or
$\{x\}\cup V\in\tau$. Hence $x\models\Box p$ or $x\models\Box q$.
This proves that $X\models(lin)$.

Suppose now $X$ is not primal. Then there exist $x\in X$ and
$U,V\in\tau$ such that $\{x\}\cup U\cup V\in\tau$, but $\{x\}\cup
U\not\in\tau$ and $\{x\}\cup V\not\in\tau$. Take a valuation such
that $P=U$ and $Q=V$. Then clearly
$x\models\Box(\Box^+p\vee\Box^+q)$. However, neither $x\models\Box
p$ nor $x\models\Box q$ is true. Indeed, if, for example,
$x\models\Box p$, then there exists an open neighborhood $W$ of $x$
such that $W\setminus\{x\}\subseteq P=U$. But then $\{x\}\cup
U=W\cup U\in\tau$, which is a contradiction. This shows that
$X\not\models(lin)$.\ep

The primal scattered spaces are very close to \emph{maximal
scattered} spaces of \cite{BekGab11}. In fact, each maximal
scattered space is primal, but there are primal spaces which are not
maximal. It follows that the logic of maximal scattered spaces is
$\ml{GL.3}$.

\section{GLP-algebras and polymodal provability logic}

A natural generalization of provability logic $\GL$ to a language
with infinitely many modalities $\la 0\ra$, $\la 1\ra$, \dots\ has
been introduced as early as in 1986 by Giorgi Japaridze
\cite{Dzh86}. He interpreted $\la 1\ra\phi$ as an arithmetical
statement expressing the \emph{$\gw$-consistency of $\phi$} over a
given g\"odelian theory $T$.\footnote{A g\"odelian theory $U$ is
$\gw$-consistent if its extension by unnested applications of the
$\gw$-rule $U':=U+\{\al{x}\phi(x): \al{n} U\vdash\phi(\ul{n})\}$ is
consistent. } Similarly, $\la 2\ra\phi$ was interpreted as the
consistency of the extension of $T+\phi$ by two nested applications
of $\gw$-rule, and so on.

While the logic of each of the individual modalities $\la n\ra$ over
Peano arithmetic was known to coincide with $\GL$ by a relatively
straightforward extension of the Solovay Theorem \cite{Boo79},
Japaridze found a complete axiomatization of the \emph{joint} logic
of the modalities $\la n\ra$, for all $n\in\gw$. This result
involved considerable technical difficulties and lead to one of the
first genuine extensions of Solovay's arithmetical fixed-point
construction. Later, Japaridze's work has been simplified and
extended by Konstantin Ignatiev \cite{Ign93} and George Boolos
\cite{Boo93a}. In particular, Ignatiev showed that $\GLP$ is
complete for more general sequences of `strong' provability
predicates in arithmetic and analyzed the variable-free fragment of
$\GLP$. Boolos included a treatment of $\GLP$ (for the somewhat
simpler case of just two modalities) in his popular book on
provability logic \cite{Boo93}.

More recently, $\GLP$ has found interesting applications in
proof-theoretic analysis of arithmetic \cite{Bek04,Bek05,AB04,Bek06}
which stimulated some further interest in the study of modal-logical
properties of $\GLP$ \cite{Car05,BJV,Ica09,Bek10}. For such
applications, the algebraic language appears to be more natural and
a different choice of the interpretation of the provability
predicates is needed. The relevant structures have been introduced
in \cite{Bek04} under the name \emph{graded provability algebras}.

Recall that an arithmetical formula is called $\Pi_n$ if it can be
obtained from a formula containing only bounded quantifiers
$\al{x\leq t}$ and $\exists{x\leq t}$ by a prefix of $n$ alternating
blocks of quantifiers starting from $\forall$. Arithmetical
$\Sigma_n$-formulas are defined dually.

Let $T$ be a g\"odelian theory. $T$ is called \emph{$n$-consistent}
if $T$ together with all true arithmetical $\Pi_n$-sentences is
consistent. (Alternatively, $T$ is $n$-consistent iff every
$\Sigma_n$-sentence provable in $T$ is true.) Let $n\nCon(T)$ denote
a natural arithmetical formula expressing the $n$-consistency of $T$
(it can be defined using the standard $\Pi_n$-definition of truth
for $\Pi_n$-sentences in arithmetic). Since we assume $T$ to be
recursively enumerable, it is easy to check that the formula
$n\nCon(T)$ itself belongs to the class $\Pi_{n+1}$.

The $n$-consistency formula induces an operator $\la n\ra_T$ acting
on the Lindenbaum--Tarski algebra $\cL_T$:
$$\la n\ra_T: [\phi]\longmapsto [n\nCon(T+\phi)]. $$
The dual \emph{$n$-provability} operators are defined by $[n]_T
x=\neg\la n\ra_T x$, for all $x\in\cL_T$. Since every true
$\Pi_n$-sentence is assumed to be an axiom for $n$-provability, we
notice that every true $\Sigma_{n+1}$-sentence must be $n$-provable.
Moreover, this latter fact is formalizable in $T$, so we obtain the
following lemma. (By the abuse of notation we denote by $[n]_T\phi$
also the arithmetical formula expressing the $n$-provability of
$\phi$ in $T$.)

\bl \label{sigman-com} For each true $\Sigma_{n+1}$-formula
$\gs(x)$,
$$T\vdash \al{x}(\gs(x)\to [n]_T\gs(\ul{x})).$$
\el

As a corollary we obtain the following basic observation probably
due to Smorynski (see~\cite{Smo85}).

\bpr For each $n\in\gw$, the structure $(\cL_T,\la n\ra_T)$ is a
Magari algebra. \epr A proof of this fact consists of verifying the
Bernays--L\"ob derivability conditions for $[n]_T$ in $T$ and of
deducing from them, in the usual way, an analog of L\"ob's theorem
for $[n]_T$.

The structure $(\cL_T,\{\la n\ra_T:n\in\gw\})$ is called the
\emph{graded provability algebra of $T$} or the \emph{GLP-algebra of
$T$}. Apart from the identities inherited from the structure of
Magari algebras for each $\la n\ra$ it satisfies the following
principles, for all $m<n$:
\begin{enumr}
\item[P1.] $\la m \ra x\leq [n]\la m \ra x$;
\item[P2.] $\la n \ra x \leq \la m \ra x$.
\end{enumr}
The validity of P1 follows from Lemma \ref{sigman-com}, because the
formula $\la m \ra_T \phi$, for any $\phi$, belongs to the class
$\Pi_{m+1}$. P2 holds, since $\la n\ra_T\phi$ asserts the
consistency of a stronger theory than $\la m\ra_T\phi$, for $m<n$.

In general, we call a \emph{GLP-algebra} a structure $(M,\{\la
n\ra:n\in\gw\})$ such that each $(M,\la n\ra)$ is a Magari algebra
and conditions P1, P2 (that are equivalent to identities) are
satisfied for all $x\in M$.

At this point it is worth noticing that P1 has an equivalent form
(modulo the other identities) that has proved to be quite useful for
various applications of GLP-algebras.

\bl \label{id-glp} Modulo the identities of Magari algebras and P2,
condition P1 is equivalent to {\em
\begin{enumr}
\item[P1$'$.] $\la n\ra y \land \la m \ra x = \la n\ra (y\land \la m\ra x)$, for all $m<n$.
\end{enumr}}
\el

\bp First, we prove P1$'$. We have $y\land \la m \ra x\leq y$, hence
$$\la n\ra (y\land \la m\ra x)\leq \la n\ra y.$$ Similarly, by P2
and transitivity, \[\la n\ra (y\land \la m\ra x)\leq \la n\ra \la m
\ra x\leq \la m\ra\la m\ra x \leq \la m\ra x.\] Hence, $\la n\ra
(y\land \la m\ra x)\leq \la n\ra y \land \la m \ra x.$ In the other
direction, by P1, $$\la n\ra y \land \la m \ra x \leq \la n\ra
y\land [n]\la m\ra x.$$ However, as in any modal algebra, we also
have $\la n \ra y \land [n] z \leq \la n\ra (y\land z)$. It follows
that
$$\la n\ra y\land [n]\la m\ra x\leq \la n \ra (y\land \la m\ra x).$$ Thus, P1$'$ is proved.

To infer P1 from P1$'$ it is sufficient to prove that $$\la m\ra x
\land \neg [n]\la m\ra x = \bot.$$ We have that $\neg [n]\la m\ra x
= \la n\ra \neg \la m \ra x$. Therefore, by P1$'$, $$\la m\ra x
\land \la n\ra \neg \la m \ra x = \la n \ra (\neg\la m\ra x \land
\la m\ra x)=\la n\ra \bot = \bot,$$ as required. \ep

An equivalent formulation of Japaridze's arithmetical completeness
theorem is that any identity of $(\cL_T,\{\la n\ra_T:n\in\gw\})$
follows from the identities of GLP-algebras \cite{Dzh86}. It is
somewhat strengthened to the current formulation in
\cite{Ign93,Bek11}.

\bt[Japaridze] Suppose $T$ is g\"odelian, $T$ contains Peano
arithmetic, and $\ch(\cL_T,\la n\ra_T)=\infty$, for each $n<\gw$.
Then, an identity holds in $(\cL_T,\{\la n\ra_T:n\in\gw\})$ iff it
holds in all GLP-algebras. \et

We note that the condition $\ch(\cL_T,\la n\ra_T)=\infty$, for each
$n\in\gw$, is equivalent to $T+n\nCon(T)$ being consistent, for each
$n\in\gw$, and is clearly necessary for the validity of Japaridze's
theorem.

The logic of all GLP-algebras can also be axiomatized as a
Hilbert-style calculus. The corresponding system $\GLP$ was
originally introduced by Japaridze. $\GLP$ is formulated in the
language of propositional logic enriched by modalities $[n]$, for
all $n\in \gw$. The axioms of $\GLP$ are those of $\GL$, formulated
for each $[n]$, as well as the two analogs of P1 and P2, for all
$m<n$:
\begin{enumr}
\item[P1.] $\la m \ra \phi\to [n]\la m \ra \phi$;
\item[P2.] $[m] \phi \to [n]\phi $.
\end{enumr}
The inference rules of $\GLP$ are modus ponens and $\phi/[n]\phi$,
for each $n\in\gw$.

For any modal formula $\phi$, $\GLP\vdash \phi$ iff the identity
$t_\phi=\top$ holds in all $\GLP$-algebras. Hence, $\GLP$ coincides
with logic of all GLP-algebras as well as with the logic of the
GLP-algebra of $T$, for any theory $T$ such that $T+n\nCon(T)$ is
consistent, for each $n<\gw$.

\section{GLP-spaces}

Topological semantics for $\GLP$ has been first considered in
\cite{BBI09}. The main difficulty in the modal-logical study of
$\GLP$ comes from the fact that it is incomplete with respect to its
relational semantics; that is, $\GLP$ is the logic of no class of
\emph{frames}. Even though a suitable class of relational
\emph{models} for which $\GLP$ is sound and complete was developed
in \cite{Bek10}, these models are sufficiently complicated and not
so easy to handle. So, it is natural to consider the generalization
of the topological semantics we have for $\GL$. As it turns out,
topological semantics provides another natural class of GLP-algebras
which is interesting in its own right, as well as because of its
analogy with the proof-theoretic GLP-algebras.

As before, we are interested in GLP-algebras of the form
$(\cP(X),\{\la n\ra:n\in\gw\})$ where $\cP(X)$ is the boolean
algebra of subsets of a given set $X$. Since each $(\cP(X),\la
n\ra)$ is a Magari algebra, the operator $\la n\ra$ is the
derivative operator with respect to some uniquely defined scattered
topology on $X$. Thus, we come to the following definition
\cite{BBI09}.

A polytopological space $(X,\{\tau_n: n\in\gw\})$ is called a
\emph{GLP-space} if the following conditions hold, for each
$n\in\gw$:

\ben \item[D0.] $(X,\tau_n)$ is a scattered topological space;
\item[D1.] For each $A\subseteq X$, $d_{\tau_n}(A)$ is $\tau_{n+1}$-open;
\item[D2.] $\tau_n\subseteq \tau_{n+1}$.
\een

We notice that the last two conditions directly correspond to
conditions P1 and P2 of GLP-algebras. By a \emph{GLP$_m$-space} we
mean a space $(X,\{\tau_n: n<m\})$ satisfying conditions D0--D2 for
the first $m$ topologies.

\bpr \label{neighbor}
\begin{enumr}
\item If $(X,\{\tau_n:n\in\gw\})$ is a GLP-space then the structure
$(\cP(X),\{d_{\tau_n}:n\in\gw\})$ is a GLP-algebra.
\item
If $(\cP(X),\{\la n\ra:n\in\gw\})$ is a GLP-algebra, then there are
uniquely defined topologies $\{\tau_n:n\in\gw\}$ on $X$ such that
$(X,\{\tau_n:n\in\gw\})$ is a GLP-space and $\la n\ra=d_{\tau_n}$,
for each $n<\gw$.
\end{enumr}
\epr

\bp (i) Suppose $(X,\{\tau_n:n\in\gw\})$ is a GLP-space. Let
$d_n:=d_{\tau_n}$ denote the corresponding derivative operators and
let $\tilde d_n$ denote its dual $\tilde d_n(A):= X\setminus d_n
(X\setminus A)$. By Theorem \ref{sim-esa} $(\cP(X),d_n)$ is a Magari
algebra, for each $n\in\gw$. Notice that $A\in\tau_n$ iff
$A\subseteq \tilde d_n A$. If $m<n$ then $d_m A\in\tau_n$, hence
$d_m A \subseteq \tilde d_n d_m A$, hence P1 holds. Since
$\tau_n\subseteq\tau_{n+1}$ we have $d_{n+1} A\subseteq d_n A$ hence
P2 holds.

 (ii) Let $(\cP(X),\{\la n\ra:n\in\gw\})$ be a GLP-algebra.
Since each of the algebras $(\cP(X),\la n\ra)$ is Magari, by Theorem
\ref{mag-scat} a scattered topology $\tau_n$ on $X$ is defined for
which $\la n\ra = d_{\tau_n}$. In fact, we have $U\in\tau_n$ iff
$U\subseteq [n]U$. We check that conditions D1 and D2 are met.

Suppose $A$ is $\tau_n$-closed, that is, $\la n\ra A\subseteq A$.
Then $\la n+1 \ra A\subseteq \la n\ra A\subseteq A$ by P2. Hence,
$U$ is $\tau_{n+1}$-closed. Thus, $\tau_n\subseteq \tau_{n+1}$.

By P1, for any set $A$ we have $\la n\ra A\subseteq [n+1]\la n\ra
A$. Hence, $d_{\tau_n}(A)=\la n\ra A\in\tau_{n+1}$. Thus,
$(X,\{\tau_n:n\in\gw\})$ is a GLP-space. \ep

\section{Derivative topology and generated GLP-spaces}

To obtain examples of GLP-spaces let us first consider the case of
two modalities. The following basic example is due to Esakia
(private communication).

\bex Consider a bitopological space $(\Omega;\tau_0,\tau_1)$ where
$\Omega$ is an ordinal, $\tau_0$ is its left topology, and $\tau_1$
is its order topology. Esakia noticed that this space is a model of
the bimodal fragment of $\GLP$, that is, in our terminology, a
GLP$_2$-space. In fact, for any $A\subseteq\Omega$ the set
$d_0(A)=(\min A,\Omega)$ is an open interval, whenever $A$ is not
empty. Hence, condition D1 holds (the other two conditions are
immediate). Esakia also noticed that such spaces can never be
complete for $\GLP$, as the linearity axiom $(.3)$ holds for $\la
0\ra$. \eex

In general, to define GLP$_n$-spaces for $n>1$, we introduce an
operation $\tau\longmapsto \tau^+$ on topologies on a given set $X$.
This operation plays a central role in the study of GLP-spaces.

Given a topological space $(X,\tau)$, let $\tau^+$ be the coarsest
topology containing $\tau$ such that each set of the form
$d_\tau(A)$, with $A\subseteq X$, is open in $\tau^+$. Thus,
$\tau^+$ is generated by $\tau$ and $\{d_\tau(A): A\subseteq X\}$.
Clearly, $\tau^+$ is the coarsest topology on $X$ such that
$(X;\tau,\tau^+)$ is a GLP$_2$-space. Sometimes we call $\tau^+$ the
\emph{derivative topology} of $(X,\tau)$.

Getting back to Esakia's example, it is easy to verify that, on any
ordinal $\Omega$, the derivative topology of the left topology
coincides with the order topology. (In fact, any open interval is an
intersection of a downset and an open upset.)

\bex Even though we are mainly interested in scattered spaces, the
derivative topology makes sense for arbitrary spaces. The reader can
check that if $\tau$ is the \emph{coarsest} topology on a set $X$
(whose open sets are just $X$ and $\emptyset$), then $\tau^+$ is the
\emph{cofinite} topology on $X$ (whose open sets are exactly the
cofinite subsets of $X$). On the other hand, if $\tau$ is the
cofinite topology then $\tau^+=\tau$. We note that the logic of the
cofinite topology on an infinite set is {\bf KD45} (see
\cite{Ste07}). \eex

For scattered spaces, $\tau^+$ is always strictly finer than $\tau$,
unless $\tau$ is discrete. We present a proof using the language of
Magari algebras.

\bpr Suppose $(X,\tau)$ is scattered. Then $d_\tau(X)$ is not open,
unless $d_\tau(X)=\emptyset$. \epr

\bp The set $d_\tau(X)$ corresponds to the element $\Diamond \top$
in the associated Magari algebra; $d_\tau(X)$ being open means
$\Diamond \top\leq \Box\Diamond \top$. By L\"ob's principle we have
$\Box\Diamond \top\leq \Box \bot = \neg\Diamond \top$. Hence,
$\Diamond \top\leq \neg\Diamond \top$, that is, $\Diamond \top =
\bot$. This means $d_\tau(X)=\emptyset$, as required. \ep

In general, we will see later that $\tau^+$ can be much finer than
$\tau$. Notice that if $\tau$ is $T_d$ then each set of the form
$d_\tau(A)$ is $\tau$-closed. Hence, it will be clopen in $\tau^+$.
Thus, $\tau^+$ is obtained by adding to $\tau$ new clopen sets. In
particular, $\tau^+$ will be zero-dimensional if so is
$\tau$.\footnote{Recall that a topological space is zero-dimensional
if it has a base of clopen sets.}

Iterating the plus operation delivers us a GLP-space. Let $(X,\tau)$
be a scattered space. Define: $\tau_0:=\tau$ and
$\tau_{n+1}:=\tau_n^+$. Then $(X,\{\tau_n:n\in\gw\})$ is a GLP-space
that will be called the GLP-space \emph{generated from $(X,\tau)$},
or simply the \emph{generated GLP-space}.

Thus, from any scattered space we can always produce a GLP-space in
a natural way. The question is whether this space will be
nontrivial, that is, whether we can guarantee that the topologies
$\tau_n$ are non-discrete.

In fact, the next observation shows that for many natural $\tau$
already the topology $\tau^+$ will be discrete. Recall that a
topological space $X$ is \emph{first-countable} if every point $x\in
X$ has a countable basis of open neighborhoods.

\bpr \label{nontriv} If $( X,\tau )$ is Hausdorff and
first-countable, then $\tau^+$ is discrete.\epr

\bp It is easy to see that if $(X,\tau)$ is first-countable and
Hausdorff, then every point $a\in d_\tau(X)$ is a (unique) limit
point of a countable sequence of points $A=\{a_n\}_{n\in\gw}$.
Hence, there is a set $A\subseteq X$ such that $d_\tau(A)=\{a\}$. By
D1 this means that $\{a\}$ is $\tau^+$-open. \ep

Thus, if $\tau$ is the order topology on a countable ordinal, then
$\tau^+$ is discrete. The same holds, for example, if $\tau$ is the
(non-scattered) topology of the real line.

We remark that the left topology $\tau$ on any countable ordinal
$>\gw$ yields an example of a first-countable space such that
$\tau^+$ is non-discrete. In the following section we will also see
that, if $\tau$ is the order topology on any ordinal $>\gw_1$, then
$\tau^+$ is non-discrete ($\gw_1$ is its least non-isolated point).
However, we do not have any topological characterization of spaces
$(X,\tau)$ such that $\tau^+$ is discrete. (See, however,
Proposition \ref{doubly-ref}, which provides a characterization in
terms of $d$-reflection.)

It is a natural question to ask what kind of separation properties
is $\tau^+$ guaranteed to have, for an arbitrary (scattered)
topology $\tau$. In fact, for $\tau^+$ we can infer a bit more
separation than for an arbitrary scattered topology. Recall that a
topological space $X$ is $T_1$ if, for any two different points
$a,b\in X$, there is a open set $U$ such that $a\in U$ and $b\notin
U$.

\bpr Let $(X,\tau)$ be any topological space. Then $(X,\tau^+)$ is
$T_1$. \epr

\bp  Let $a,b\in X$, $a\neq b$. We must show that there is an open
set $A$ such that $a\in A$ and $b\notin A$. Consider the set
$B:=d_\tau(\{b\})$, which is open in $\tau^+$. If $a\in B$ then we
are done, because $b\notin d_\tau(\{b\})=B$. Otherwise, if $a\notin
B$, neither does $a$ belong to the closure of $\{b\}$ (which is
simply $\{b\}\cup d_\tau(\{b\})$). It follows that the complement of
$\{b\}\cup B$ is the required open set. \ep

The following example shows that, in general, $\tau^+$ need not
always be Hausdorff.

\bex \label{t-one} Let $(X,\prec)$ be a strict partial ordering on
$X:=\gw\cup\{a,b\}$ where $\gw$ is taken with its natural order, $a$
and $b$ are $\prec$-incomparable, and $n\prec a,b$, for all
$n\in\gw$. Let $\tau$ be the left topology on $(X,\prec)$. Since
$\prec$ is well-founded, $\tau$ is scattered.

Notice that $d_\tau(A)=\{x\in X:\exists y\in A\ y\prec x\}.$ Hence,
if $A$ intersects $\gw$, then $d_\tau(A)$ contains an end-segment of
$\gw$. Otherwise, $d_\tau(A)=\emptyset$. It follows that a base of
open neighborhoods of $a$ in $\tau^+$ consists of sets of the form
$I\cup\{a\}$ where $I$ is an end-segment of $\gw$. Similarly, sets
of the form $I\cup \{b\}$ are a base of open neighborhoods of $b$.
But any two such sets have a non-empty intersection. \eex

\section{$d$-reflection} \label{d-ref}

In the next section we are going to describe in some detail the
GLP-space generated from the left topology on the ordinals.
Strikingly, we will see that it naturally leads to some of the
central notions of combinatorial set theory, such as Mahlo operation
and stationary reflection. In fact, part of our analysis can be
easily stated using the language of modal logic for arbitrary
generated GLP-spaces. In this section we provide a necessary setup
and characterize the topologies of a generated GLP-space in terms of
what we call \emph{$d$-reflection}. \footnote{Curiously, the reader
may notice that the notion of \emph{reflection principle} as used in
provability logic and formal arithmetic matches very nicely the
notions such as \emph{stationary reflection} in set theory. (As far
as we know, the two terms have evolved completely independently from
one another.)}

Throughout this section we will assume that $(X,\tau)$ is a $T_d$
space and let $d=d_\tau$. Thus, $(\cP(X),d)$ is a K4-algebra. \bd A
point $a\in X$ is called \emph{$d$-reflexive} if $a\in dX$ and, for
each $A\subseteq X$,
\[a\in dA \Imp a\in d(dA).\]
In modal logic terms this means that the following formula is valid
at $a\in X$, for any evaluation of the variable $p$:
\[X,a\models \Diamond \top \land (\Diamond p\to \Diamond\Diamond p).\]
Similarly, a point $a\in X$ is called \emph{$m$-fold $d$-reflexive}
if $a\in dX$ and, for each $A_1,\dots,A_m\subseteq X$,
\[a\in dA_1\cap\dots\cap d A_m \Imp a\in d(dA_1\cap\dots\cap d A_m).\]
$2$-fold $d$-reflexive points will also be called \emph{doubly
$d$-reflexive} points. Modal logically, $a\in X$ is doubly
$d$-reflexive iff
\[X,a\models \Diamond \top \land (\Diamond p\land \Diamond q\to \Diamond(\Diamond p\land \Diamond q)).\]
\ed

\bl \label{m-fold} Each doubly $d$-reflexive point $x\in X$ is
$m$-fold $d$-reflexive, for any finite $m$. \el

\bp The argument goes by induction on $m\geq 2$. Suppose $x\in
dA_1\cap \dots \cap dA_{m+1}$, then $x\in dA_1\cap \dots \cap
dA_{m}$ and $x\in dA_{m+1}$. By induction hypothesis, $$x\in
d(dA_1\cap \dots \cap dA_{m})$$ and by 2-fold reflection $$x\in
d(d(dA_1\cap \dots \cap dA_m)\cap dA_{m+1}).$$ However, by $T_d$
property $$d(dA_1\cap \dots \cap dA_m)\subseteq dA_1\cap \dots \cap
dA_m,$$ hence
$$x\in d(dA_1\cap \dots \cap dA_m\cap dA_{m+1}),$$ as required.\ep

\bpr \label{doubly-ref} Let $(X,\tau)$ be a $T_d$-space. A point
$x\in X$ is doubly $d$-reflexive iff $x$ is a limit point of
$(X,\tau^+)$.\epr

\bp For the (if) direction, we give an argument in the algebraic
format. In fact, it is sufficient to show the following inequality
in the algebra of $(X,\tau)$, for any elements $p,q\subseteq X$:
$$\la 1\ra \top \land\la 0\ra p \land \la 0 \ra q \leq \la 0\ra(\la 0\ra p \land \la 0 \ra
q).$$

Notice that by Proposition \ref{id-glp}, $$\la 1\ra \top \land\la
0\ra p = \la 1\ra (\top \land \la 0\ra p) = \la 1\ra\la 0\ra p.$$
Hence, we obtain using P1$'$ once again: $$\la 1\ra \top \land\la
0\ra p \land \la 0 \ra q=\la 1\ra\la 0\ra p\land \la 0\ra q = \la
1\ra(\la 0\ra p\land \la 0\ra q).$$ The latter formula can be
weakened to $\la 0\ra(\la 0\ra p\land \la 0\ra q)$ by P2, as
required.

\medskip
For the (only if) direction, it is sufficient to show that each
doubly $d$-reflexive point of $(X,\tau)$ is a limit point of
$\tau^+$. Suppose $x$ is doubly $d$-reflexive. By Lemma
\ref{m-fold}, $x$ is $m$-fold $d$-reflexive.

Any basic open subset of $\tau^+$ has the form
$$U:=A_0\cap d A_1\cap \dots \cap d A_m,$$
where $A_0\in\tau$. Assume $x\in U$, we have to find a point $y\neq
x$ such that $y\in U$.

Since $x\in d A_1\cap \dots \cap d A_m,$ by $m$-fold $d$-reflexivity
we obtain $x\in d(d A_1\cap \dots \cap d A_m).$ Since $A_0$ is an
open neighborhood of $x$, there is a $y\in A_0$ such that $y\neq x$
and $y\in d A_1\cap \dots \cap d A_m.$ Hence $y\in U$ and $y\neq x$,
as required. \ep

Let $d^+$ denote the derivative operator associated with $\tau^+$.
We obtain the following characterization of derived topology in
terms of neighborhoods.

\bpr \label{neighb} A subset $U\subseteq X$ contains a
$\tau^+$-neighborhood of $x\in X$ iff one of the following two cases
holds:

\begin{enumr}
\item $x$ is not doubly $d$-reflexive and $x\in U$;
\item $x$ is doubly $d$-reflexive and there is an $A\in \tau$ and a $B$ such that $x\in A\cap dB\subseteq U$.
\end{enumr}
\epr

\bp Since (i) ensures that $x$ is $\tau^+$-isolated, each condition
is clearly sufficient for $U$ to contain a $\tau^+$-neighborhood of
$x$. To prove the converse, assume that $U$ contains a
$\tau^+$-neighborhood of $x$. This means $$x\in A\cap dA_1\cap \dots
\cap dA_m \subseteq U,$$ for some $A,A_1, \dots, A_m$ with $A\in
\tau$.

If $x$ is $\tau^+$-isolated, condition (i) holds. Otherwise, $x\in
d^+ X$. Let $B:=dA_1\cap\dots \cap dA_m$. Since $B$ is closed in
$\tau$ we have $dB\subseteq B$, hence $A\cap dB\subseteq U$. It
remains us to show that $x\in A\cap dB$. By Lemma \ref{id-glp}, $B
\cap d^+ X = d^+ B \subseteq dB.$ Hence, $x\in A\cap B \cap d^+
X\subseteq A\cap dB$. \ep

\brem Since in (ii) of the previous lemma $A$ is open, we have
$A\cap dB= A\cap d(A\cap B)$, for any $B$. Hence, we can assume
$B\subseteq A$. \erem

\bc \label{d-plus} Let $(X,\tau)$ be a $T_d$-space. Then, for all
$x\in X$ and $A\subseteq X$, $x\in d^+ A$ iff  the following two
conditions hold:
\begin{enumr}
\item $x$ is doubly $d$-reflexive;
\item For all $B\subseteq X$, $x\in dB \Imp x\in d(A\cap dB)$.
\end{enumr}
\ec The second condition is similar to $x$ being $d$-reflexive,
however $dB$ is now required to reflect from $x$ to a point of $A$.

\bp The fact that (i) and (ii) are necessary is proved using
Proposition \ref{doubly-ref} and the inequality
$$d^+ A \cap d B = d^+ (A \cap dB)\subseteq d(A\cap dB).$$
We prove that (i) and (ii) are sufficient. Assume $x\in U\in
\tau^+$. By Proposition \ref{neighb} we can assume that $U$ has the
form $V\cap dB$ where $V\in\tau$. By (ii), from $x\in dB$ we obtain
$x\in d(A\cap dB)$. Hence, there is a $y\neq x$ such that $y\in V$
and $y\in A\cap dB$. It follows that $y\in A$ and $y\in V\cap dB=U$.
\ep

\section{The ordinal GLP-space} \label{ordinal}

Here we discuss the GLP-space generated from the left topology on
the ordinals, that is, the GLP-space $(\Omega;\{\tau_n:n\in\gw\})$
with $\tau_0$ the left topology on $\Omega$ and
$\tau_{n+1}=\tau_n^+$, for each $n\in\gw$. It is convenient to think
of $\Omega$ as the class of all ordinals, even though some readers
might feel safer with $\Omega$ being a large fixed ordinal. Our
basic findings are summarized in the following table, to which we
provide extended comments below.

The rows of the table correspond to topologies $\tau_n$. The first
column contains the name of the topology (the first two are
standard, the third one is introduced in \cite{BBI09}, the fourth
one is introduced here). The second column indicates the first limit
point of $\tau_n$, which is denoted $\theta_n$. The last column
describes the derivative operator associated with $\tau_n$. We note
that $\theta_3$ is a large cardinal which is sometimes referred to
as \emph{the first cardinal reflecting for pairs of stationary sets}
(see below), but we know no special notation for this cardinal.

\bigskip
\begin{tabular}{r|c|c|l}
 & name & $\theta_n$ &\ \  $d_n(A)$  \\ \hline
$\tau_0$ \ & left & 1 & $\{\ga: A\cap\ga\neq\emptyset\}$  \\ \hline
$\tau_1$ \ & order & $\gw$ & $\{\ga\in Lim:\text{$A\cap\ga$ is
unbounded in $\ga$}\}$ \\ \hline
$\tau_2\ $ & club & $\gw_1$ &
$\{\ga:\text{$\cf(\ga)>\gw$ and $A\cap\ga$ is stationary in
$\ga$}\}$ \\ \hline
$\tau_3\ $ &\ \ \ Mahlo\ \ \ & $\ \ \theta_3\ \ $ & \dots\ \dots
\end{tabular}

\bigskip

We have already seen that the derivative topology of the left
topology is exactly the order topology. Therefore, basic facts
related to the first two rows of the table are rather clear. We turn
to the next topology $\tau_2$.

\subsection{Club topology}

To characterize $\tau_2$ we apply Proposition \ref{neighb}, hence it
is useful to see what corresponds to the notion of doubly
$d$-reflexive point of the interval topology.

Recall that the \emph{cofinality} $\cf(\ga)$ of a limit ordinal
$\ga$ is the least order type of a cofinal subset of $\ga$;
$\cf(\ga):=0$ if $\ga\notin\Lim$. An ordinal $\ga$ is \emph{regular}
if $\cf(\ga)=\ga$.

\bl \label{cofin} For any ordinal $\ga$, $\ga$ is $d_1$-reflexive
iff $\ga$ is doubly $d_1$-reflexive iff $\cf(\ga)>\gw$. \el

\bp $d_1$-reflexivity of $\ga$ implies $\cf(\ga)>\gw$. In fact,
$d_1$-reflexivity of $\ga$ means that $\ga\in \Lim$ and, for all
$A\subseteq \ga$, if $A$ is cofinal in $\ga$, then $d_1(A)$ is
cofinal in $\ga$. If $\cf(\ga)=\gw$ then there is an increasing
sequence $(\ga_n)_{n\in\gw}$ such that $\sup\{\ga_n:n\in\gw\}=\ga$.
Then, for $A:=\{\ga_n:n\in\gw\}$ we obviously have $d_1(A)=\{\ga\}$,
hence $A$ violates the reflexivity property.

Now we show that $\cf(\ga)>\gw$ implies $\ga$ is doubly
$d_1$-reflexive. Suppose $\cf(\ga)>\gw$ and $A,B\subseteq \ga$ are
both cofinal in $\ga$. We show that $d_1A\cap d_1B$ is cofinal in
$\ga$. Assume $\gb<\ga$. Using the cofinality of $A,B$ we can
construct an increasing sequence $(\gy_n)_{n\in\gw}$ above $\gb$
such that $\gy_n\in A$, for even $n$, and $\gy_n\in B$ for odd $n$.
Let $\gy:=\sup \{\gy_n:n<\gw\}$. Obviously, both $A$ and $B$ are
cofinal in $\gy$ whence $\gy\in d_1A\cap d_1B$. Since $\cf(\ga)>\gw$
and $\cf(\gy)=\gw$, we have $\gy<\ga$. \ep

\bc Limit points of $\tau_2$ are exactly the ordinals of uncountable
cofinality. \ec

It turns out that topology $\tau_2$ is strongly related to the
well-known concept of \emph{club filter} in set theory.

Let $\ga$ be a limit ordinal. A subset $C\subseteq \ga$ is called a
\emph{club} in $\ga$ if $C$ is closed in the order topology of $\ga$
and unbounded in $\ga$. The filter on $\ga$ generated by all clubs
in $\ga$ is called the \emph{club filter}.


\bpr \label{club} Assume $\cf(\ga)>\gw$. The following statements
are equivalent:
\begin{enumr}
\item $U$ contains a $\tau_2$-neighborhood of $\ga$;
\item There is a $B\subseteq\ga$ such that $\ga\in d_1B\subseteq U$;
\item $\ga\in U$ and $U$ contains a club in $\ga$;
\item $\ga\in U$ and $U\cap \ga$ belongs to the club filter on $\ga$.
\end{enumr}
\epr

\bp  Statement (ii) implies (iii), since $\ga\cap d_1 B$ is a club
in $\ga$, whenever $\ga\in d_1B$. Statement (iii) implies (iv), for
obvious reasons.

Statement (iv) implies (i). If $C$ is a club in $\ga$, then
$C\cup\{\ga\}$ contains a $\tau_2$-neighborhood $d_1 C$ of $\ga$.
Indeed, $d_1 C$ is $\tau_2$-open, contains $\ga$, and $d_1
C\subseteq C\cup \{\ga\}$ since $C$ is $\tau_1$-closed in $\ga$.

Statement (i) implies (ii). Assume $U$ contains a
$\tau_2$-neighborhood of $\ga$. Since $\cf(\ga)>\gw$, by Lemma
\ref{cofin} and Proposition \ref{neighb} there is an $A\in \tau_1$
and a $B_1$ such that $\ga\in A\cap d_1  B_1\subseteq U$. Since $A$
is a $\tau_1$-neighborhood of $\ga$, by Proposition \ref{neighb}
there are $A_0\in\tau_0$ and a $B_0$ such that $\ga\in A_0\cap d_0
B_0$. Since $\tau_0$ is the left topology, we can assume that $A_0$
is the minimal $\tau_0$-neighborhood $[0,\ga]$ of $\ga$. Besides, we
have $\ga\in d_0B_0 \cap d_1 B_1= d_1(B_1 \cap d_0B_0)\subseteq U$.
Since $[0,\ga]$ is $\tau_1$-clopen, $d_1(C\cap \ga)=[0,\ga]\cap d_1
C$ for any $C$, so we can take $B_1\cap d_0B_0\cap\ga$ for $B$.\ep

\bd The \emph{club topology} on $\Omega$ is the unique topology such
that \bi \item If $\cf(\ga)\leq \gw$ then $\ga$ is an isolated
point;
\item If $\cf(\ga)>\gw$, then, for any $U\subseteq\Omega$, $U$ contains a neighborhood of $\ga$ iff $\ga\in U$ and $U$ contains a club in $\ga$.
\ei \ed

As an immediate corollary of Proposition \ref{club} we obtain \bt
\label{tau-two} $\tau_2$ coincides with the club topology. \et

We remark that the above theorem saves us the little work of
verifying that the club topology is, indeed, a topology.

The derivative operation for the club topology is also well-known in
set theory. Recall the following definition for $\cf(\ga)>\gw$.

A subset $A\subseteq \ga$ is called \emph{stationary in $\ga$} if
$A$ intersects every club in $\ga$. Observe that this happens
exactly when $\ga$ is a limit point of $A$ in $\tau_2$, so
$$d_2(A)=\{\ga: \text{ $\cf(\ga)>\gw$ and $A\cap\ga$ is stationary in $\ga$}\}.$$
The map $d_2$ is usually called the \emph{Mahlo operation} (see
\cite{Jech} where $d_2$ is denoted Tr). Its main significance is
associated with the notion of Mahlo cardinal, one of the basic
examples of large cardinals in set theory. Let $\Reg$ denote the
class of regular cardinals; the ordinals in $d_2(\Reg)$ are called
\emph{weakly Mahlo cardinals}. Their existence implies the
consistency of $\ZFC$, as well as the consistency of $\ZFC$ together
with the assertion `inaccessible cardinals exist.'

Now we turn to topology $\tau_3$.

\subsection{Stationary reflection and Mahlo topology}

Since the open sets of $\tau_3$ are generated by the Mahlo
operation, we call $\tau_3$ \emph{Mahlo topology}. It turns out to
be intrinsically connected with the concept of stationary
reflection, an extensively studied concept in set theory.

We adopt the following terminology. An ordinal $\lambda$ is called
\emph{reflecting} if $\cf(\lambda)>\gw$ and, whenever $A$ is
stationary in $\lambda$, there is an $\ga<\lambda$ such that $A\cap
\ga$ is stationary in $\ga$. Similarly, $\lambda$ is \emph{doubly
reflecting} if $\cf(\lambda)>\gw$ and whenever $A,B$ are stationary
in $\lambda$ there is an $\ga<\lambda$ such that both $A\cap \ga$
and $B\cap \ga$ are stationary in $\ga$.

Mekler and Shelah's notion of \emph{reflection cardinal}
\cite{MekShe} is somewhat more general than the one given here,
however it has the same consistency strength. Reflection for pairs
of stationary sets has been introduced by Magidor \cite{Magi82}.
Since $d_2$ coincides with the Mahlo operation, we immediately
obtain the following statement.

\bpr
\begin{enumr}
\item $\lambda$ is reflecting iff $\lambda$ is $d_2$-reflexive;
\item $\lambda$ is doubly reflecting iff $\lambda$ is doubly $d_2$-reflexive;
\item $\lambda$ is a non-isolated point in $\tau_3$ iff $\lambda$ is doubly reflecting.
\end{enumr}
\epr

Together with the next proposition this yields a characterization of
Mahlo topology in terms of neighborhoods.

\bpr Suppose $\gl$ is doubly reflecting. For any subset $U\subseteq
\Omega$, the following conditions are equivalent:
\begin{enumr}
\item $U$ contains a $\tau_3$-neighborhood of $\gl$;
\item $\gl\in U$ and there is a $B\subseteq \gl$ such that $\gl\in d_2 B\subseteq U$;
\item $\gl\in U$ and there is a $\tau_2$-closed (in the relative topology of $\gl$) stationary $C\subseteq \gl$ such that $C\subseteq U$.
\end{enumr}
\epr

Notice that the notion of $\tau_2$-closed stationary $C$ in (iii) is
the analog of the notion of club for the $\tau_2$-topology.

\bp Condition (ii) implies (iii). Since $\gl$ is reflecting, if
$\gl\in d_2 B$ then $\gl\in d_2 d_2 B$, that is, $\gl\cap d_2 B$ is
stationary in $\gl$. So we can take $C:= \gl\cap d_2 B$.

Condition (iii) implies (ii). If $C$ is $\tau_2$-closed and
stationary in $\gl$, then $d_2 C\subseteq C\cup\{\gl\}\subseteq U$
and $\gl\in d_2 C$. Thus, $\gl\cap d_2 C$ can be taken for $B$.

Condition (ii) implies (i). If (ii) holds, $U$ contains a subset of
the form $d_2 B$. The latter is $\tau_3$-open and contains $\gl$,
thus, is a neighborhood of $\gl$.

For the converse direction, we note that by Proposition \ref{neighb}
$U$ contains a subset of the form $A\cap d_2 B$ where $A\in \tau_2$,
$B\subseteq A$ and $\gl\in A\cap d_2 B$. Since $A$ is a
$\tau_2$-neighborhood of $\gl$, by Theorem \ref{tau-two} there is a
set $B_1$ such that $\gl\in [0,\gl]\cap d_1 B_1\subseteq A$. Then,
$$\gl\in [0,\gl]\cap d_1 B_1\cap d_2 B = [0,\gl]\cap d_2 (B\cap
d_1B_1).$$ Since $[0,\gl]$ is clopen, we obtain $\gl\in d_2 C$ with
$C:= B\cap d_1 B_1\cap\gl$. \ep

Reflecting and doubly reflecting cardinals are large cardinals in
the sense that their existence implies consistency of $\ZFC$. They
have been studied by Mekler and Shelah \cite{MekShe} and Magidor
\cite{Magi82} who investigated their consistency strength and
related them to some other well-known large cardinals. By a result
of Magidor, the existence of a doubly reflecting cardinal is
equiconsistent with the existence of a \emph{weakly compact
cardinal}.\footnote{Weakly compact cardinals are the same as
$\Pi_1^1$-indescribable cardinals, see below.} More precisely, the
following proposition holds.

\bpr \begin{enumr} \item If $\lambda$ is weakly compact then
$\lambda$ is doubly reflecting. \item \emph{(Magidor)} If $\lambda$
is doubly reflecting then $\lambda$ is weakly compact in $L$.
\end{enumr} \epr

Here, the first item is well-known and easy. Magidor originally
proved the analog of the second item for $\gl=\aleph_2$ and
stationary sets of ordinals of countable cofinality in $\aleph_2$.
However, it has been remarked by Mekler and Shelah \cite{MekShe}
that essentially the same proof yields the stated
claim.\footnote{The first author thanks J. Cummings for clarifying
this.}

\bc Assertion ``$\tau_3$ is non-discrete'' is equiconsistent with
the existence of a weakly compact cardinal.\ec

\bc \label{discr} If $\ZFC$ is consistent then it is consistent with
$\ZFC$ that $\tau_3$ is discrete and hence that $\GLP_3$ is
incomplete w.r.t.\ any ordinal space.\ec

Recall that $\theta_n$ denotes the first non-isolated point of
$\tau_n$ (in the space of all ordinals). We have: $\theta_0=1$,
$\theta_1=\gw$, $\theta_2=\gw_1$, $\theta_3$ is the first doubly
reflecting cardinal.

$\ZFC$ does not know much about the location of $\theta_3$, however
the following facts are interesting.

\bi \item $\theta_3$ is regular, but not a successor of a regular
cardinal;

\item While weakly compact cardinals are non-isolated,
$\theta_3$ need not be weakly compact: If infinitely many
supercompact cardinals exist, then there is a model where
$\aleph_{\gw+1}$ is doubly reflecting (Magidor \cite{Magi82});
\item If $\theta_3$ is a successor of a singular cardinal, then some very
strong large cardinal hypothesis must be consistent (Woodin
cardinals). \ei

\subsection{Further topologies}

Further topologies of the ordinal GLP-space do not seem to have
prominently occured in set-theoretic work. They yield some large
cardinal notions, for the statement that $\tau_n$ is non-discrete
(or, equivalently, the statement that $\theta_n$ exists) implies the
existence of a doubly reflecting cardinal, for any $n>2$. We do not
know whether cardinals $\theta_n$ coincide with any of the standard
large cardinal notions.

Here we give a sufficient condition for the topology $\tau_{n+2}$ to
be non-discrete. We show that, if there exists a
\emph{$\Pi_n^1$-indescribable cardinal}, then $\tau_{n+2}$ is
non-discrete.

Let $Q$ be a class of second order formulas over the standard first
order set-theoretic language enriched by a unary predicate $R$. We
assume  $Q$ to contain at least the class of all first order
formulas (denoted $\Pi^1_0$). We shall consider standard models of
that language of the form $(V_\ga,\in,R)$, where $\ga$ is an
ordinal, $V_\ga$ is the $\ga$-th class in the cumulative hierarchy,
and $R$ is a subset of $V_\ga$.

We would like to give a definition of $Q$-indescribable cardinals in
topological terms. They can then be defined as follows.

\bd For any sentence $\phi\in Q$ and any $R\subseteq V_\kappa$, let
$U_\kappa(\phi,R)$ denote the set $$\{\ga\leq\kappa:
(V_\ga,\in,R\cap V_\ga)\models \phi\}.$$

The \emph{$Q$-describable topology $\tau_Q$} on $\Omega$ is
generated by a subbase consisting of sets $U_\kappa(\phi,R)$, for
all $\kappa\in\Omega$, $\phi\in Q$, and all $R\subseteq V_\kappa$.
\ed

As an exercise, the reader can check that the intervals $(\ga,\gb]$
are open in any $\tau_Q$ (consider $R=\{\ga\}$ and
$\phi=\exi{x}(x\in R)$). The main strength of the $Q$-describable
topology, however, comes from the fact that a second order variable
$R$ is allowed to occur in $\phi$. So, a lot of subsets of $\Omega$
that can be `described' in this way are open in $\tau_Q$.

Let $d_Q$ denote the derivative operator for $\tau_Q$. An ordinal
$\kappa$ is called \emph{$Q$-indescribable} if it is a limit point
of $\tau_Q$. In other words, $\kappa$ is $Q$-indescribable iff
$\kappa\in d_Q(\On)$ iff $\kappa \in d_Q(\kappa)$.

It is not difficult to show that, whenever $Q$ is any of the classes
$\Pi_n^1$, the sets $U_\kappa(\phi,R)$ actually form a base for
$\tau_Q$. Hence, our definition of $\Pi_n^1$-indescribable cardinals
is equivalent to the standard one given in \cite{Kan09}: $\kappa$ is
$Q$-indescribable iff,
 for all $R\subseteq
V_\kappa$ and all sentences $\phi\in Q$,
$$(V_\kappa,\in,R)\models \phi\ \Imp \ \exists \ga<\kappa\: (V_\ga,\in,R\cap
V_\ga)\models \phi.$$

A \emph{weakly compact cardinal} can be defined as the
$\Pi_1^1$-indescribable one. The well-known fact that weakly compact
cardinals are doubly reflecting can be somewhat more generally
stated as follows: the Mahlo topology $\tau_3$ is contained in
$\tau_{\Pi_1^1}$. We omit a short proof, because we are going to
prove a more general proposition suggested to the first author by
Philipp Schlicht (see \cite{Bek09a}).

\bpr For any $n\geq 0$, $\tau_{n+2}$ is contained in
$\tau_{\Pi_n^1}$. \epr

\bp We shall show that, for each $n$, there is a $\Pi_n^1$-formula
$\phi_{n+1}(R)$ such that
$$\kappa\in  d_{n+1}(A) \iff (V_\kappa,\in, A\cap \kappa)\models\phi_{n+1}(R). \eqno (**)$$
This implies that, for each $\kappa\in  d_{n+1}(A)$, the set
$U_\kappa(\phi_{n+1},A\cap\kappa)$ is a $\tau_{\Pi_n^1}$-open subset
of $ d_{n+1}(A)$ containing $\kappa$. Hence, each $ d_{n+1}(A)$ is
$\tau_{\Pi_n^1}$-open. Since $\tau_{n+2}$ is generated over
$\tau_{n+1}$ by the open sets of the form $ d_{n+1}(A)$ for various
$A$, we have $\tau_{n+2}\subseteq\tau_{\Pi_n^1}$.

We prove $(**)$ by induction on $n$. For $n=0$, notice that
$\kappa\in d_1(A)$ iff ($\kappa\in\Lim$ and $A\cap\kappa$ is
unbounded in $\kappa$) iff
$$(V_\kappa,\in,A\cap \kappa)\models \al{\ga}\exists \gb\: (R(\gb) \land \ga<\gb).$$

For the induction step recall that, by Corollary \ref{d-plus},
$\kappa\in  d_{n+1}(A)$ iff
\begin{enumr}
\item $\kappa$ is doubly $ d_n$-reflexive; \label{one}
\item $\forall Y\subseteq \kappa\:(\kappa\in d_n(Y)\to \exists
\ga<\kappa\:(\ga\in A \land \ga\in  d_n(Y))$.\label{two}
\end{enumr}
By the induction hypothesis, for some $\phi_n(R)\in\Pi_{n-1}^1$, we
have
$$\alpha\in  d_n(A) \iff (V_\alpha,\in, A\cap \ga)\models\phi_n(R).$$
Hence, part (ii) is equivalent to
$$(V_\kappa,\in,A\cap\kappa)\models \forall Y\subseteq\On\,(\phi_n(Y)\to \exists
\ga\:(R(\ga) \land \phi_n^{V_\ga}(Y\cap \ga))).$$ Here,
$\phi^{V_\ga}$ means the relativization of all quantifiers in $\phi$
to $V_\ga$. We notice that $V_\ga$ is first order definable, hence
the complexity of $\phi_n^{V_\ga}$ remains in the class
$\Pi_{n-1}^1$. So, the resulting formula is $\Pi_n^1$.

To treat part (i) we recall that $\kappa$ is doubly $ d_n$-reflexive
iff $\kappa\in d_n(\On)$ and
$$\forall Y_1,Y_2\subseteq\kappa\,(\kappa\in d_n(Y_1)\cap d_n(Y_2) \to
\exists \ga<\kappa\:\ga\in  d_n(Y_1)\cap d_n(Y_2)).$$ Similarly to
the above, using the induction hypothesis this can be rewritten as a
$\Pi_n^1$-formula. \ep

\bc If there is a $\Pi_n^1$-indescribable cardinal, then
$\tau_{n+2}$ has a non-isolated point. \ec

\bc If there is a cardinal which is $\Pi_n^1$-indescribable, for
each $n$, then all $\tau_n$ are non-discrete.  \ec

We know that $\theta_3$ need not be weakly compact in some models of
$\ZFC$. Hence, the condition of the existence of
$\Pi_n^1$-indescribable cardinals is not a necessary one for the
nontriviality of the topologies $\tau_{n+2}$, in general. However,
the two statements are equiconsistent, as has been recently shown by Bagaria, Magidor, and Sakai in \cite{BaMaSa12}.

\section{Topological completeness results}

As in the case of the unimodal language (cf Section \ref{topcompl}),
one can ask two basic questions: \bi
\item Is $\GLP$ complete w.r.t.\ the class of all GLP-spaces?
\item Is $\GLP$ complete w.r.t.\ some fixed natural GLP-spaces?
\ei

In the unimodal case, both questions received positive answers due
to Esakia and Abashidze--Blass, respectively. Now the situation is
more complicated.

\subsection{Topological completeness of $\GLP$}

The first question was initially studied by Beklemishev,
Bezhanishvili and Icard in \cite{BBI09}. In this paper only some
partial results were obtained. It was proved that the bimodal system
$\GLP_2$ is complete w.r.t.\ GLP$_2$-spaces of the form
$(X,\tau,\tau^+)$, where $X$ is a well-founded partial ordering and
$\tau$ is its left topology. A proof of this result was based on the
Kripke model techniques coming from \cite{Bek10}.

Already at that time it was clear that these techniques cannot be
immediately generalized to GLP$_3$-spaces, since the third topology
$\tau^{++}$ on such orderings is sufficiently similar to the club
topology. From the results of Andreas Blass described below it was
known that some stronger set-theoretic assumptions would be needed
to prove the completeness w.r.t.\ such topologies. Moreover, without
any large cardinal assumptions it was not even known whether a
GLP-space with a non-discrete third topology could exist at all.

First examples of GLP-spaces in which all topologies are
non-discrete are constructed in \cite{BekGab11}. In the same paper
also the much stronger fact of topological completeness of $\GLP$
w.r.t.\ the class of all (countable, Hausdorff) GLP-spaces is
established.

\bt
\begin{enumr}
\item $\Log(\cC)=\GLP$, where $\cC$ is the class of all GLP-spaces.
\item There is a countable Hausdorff GLP-space $X$ such that $\Log(X)=\GLP$.
\end{enumr}
\et

In fact, $X$ is the ordinal $\ge_0$ equipped with a sequence of
topologies refining the order topology. However, these topologies
cannot be first-countable and are, in fact, defined using
non-constructive methods such as Zorn's lemma. In this sense, it is
not an example of a \emph{natural} GLP-space.

A proof of this theorem introduces the techniques of maximal
rank-preserving and limit-maximal extensions of scattered spaces. It
falls outside the present survey (see \cite{BekGab11}).

\subsection{Completeness w.r.t.\ the ordinal GLP-space}

The question whether $\GLP$ is complete w.r.t.\ some natural
GLP-space is still open. Some partial results concerning the
GLP-space generated from the order topology on the ordinals are
described below. Here, we call this space the \emph{ordinal
GLP-space}. (The space described in the previous section is clearly
not an exact model of $\GLP$, as the left topology validates the
linearity axiom.)

As we know from Corollary \ref{discr}, it is consistent with $\ZFC$
that the Mahlo topology is discrete. Hence, it is consistent that
$\GLP$ is incomplete w.r.t.\ the ordinal GLP-space. However, is it
consistent that $\GLP$ be complete? To this question we do not know
a full answer. A pioneering work has been done by Andreas Blass
\cite{Bla90} who studied the question of completeness of the
G\"odel--L\"ob logic $\GL$ w.r.t.\ a semantics equivalent to the
topological interpretation w.r.t.\ the \emph{club topology}
$\tau_2$. He used the language of filters rather than that of
topological spaces, as is more common in set theory.

\bt[Blass] \begin{enumr} \item If $V=L$ and $\Omega\geq\aleph_\gw$,
then $\GL$ is complete w.r.t.\ $(\Omega,\tau_2)$.

\item If there is a weakly Mahlo cardinal, there is a model of $\ZFC$ in which $\GL$ is incomplete w.r.t.\ $(\Omega,\tau_2)$, for any $\Omega$.
\end{enumr}
\et

A corollary of (i) is that the statement ``$\GL$ is complete w.r.t.\
$\tau_2$'' is consistent with $\ZFC$ (provided $\ZFC$ is
consistent). In fact, instead of $V=L$ Blass used the so-called
\emph{square principle} for all $\aleph_n$, $n<\gw$, which holds in
$L$ by the results of Ronald Jensen. A proof of (i) is based on an
interesting combinatorial construction using the techniques of
splitting stationary sets.

A proof of (ii) is much easier. It uses a model of Harrington and
Shelah in which $\aleph_2$ is reflecting for stationary sets of
ordinals of countable cofinality \cite{HarShe}. Assuming Mahlo
cardinals exist, they have shown that the following statement holds
in some model of $\ZFC$:

\begin{quote} If $S$ is a stationary subset of $\aleph_2$ such that $\al{\ga\in S} \cf(\ga)=\gw$, then there is a $\gb<\ga$ (of cofinality $\gw_1$) such that $S\cap \gb$ is stationary in $\gb$.
\end{quote}

In fact, this statement can be expressed in the language of modal
logic. First, we remark that this principle implies its
generalization to all ordinals $\gl$ of cofinality $\aleph_2$
(consider an increasing continuous function mapping $\aleph_2$ to a
club in $\gl$). Second, we remark that, for the club topology, the
formula $\Diamond^n \top$ represents the class ordinals of
cofinality at least $\aleph_n$. This is a straightforward
generalization of Lemma \ref{cofin}. Thus, the formula $\Box^3\bot
\land \Diamond^2\top$ represents the subclass of $\Omega$ consisting
of ordinals of cofinality $\gw_2$.

Hence, the above reflection principle  amounts to the validity of
the following modal formula:
$$\Box^3\bot \land \Diamond^2\top\land \Diamond(p\land\Box\bot) \to \Diamond^2(p\land \Box\bot). \eqno (*)$$
In fact, if the left hand side is valid in $\gl$, then
$\cf(\gl)=\gw_2$ and the interpretation of $p\land \Box \bot$ is a
set $S$ consisting of ordinals of countable cofinality such that
$S\cap\gl$ is stationary in $\gl$. The right hand side just states
that this set reflects.  Thus, formula $(*)$ is valid in
$(\Omega,\tau_2)$, for any $\Omega$. Since this formula is clearly
not provable in $\GL$, the topological completeness fails for
$(\Omega,\tau_2)$.

Thus, Blass managed to give an exact consistency strength of the
statement ``\,$\GL$ is incomplete w.r.t.\ $\tau_2$''.

\bc ``\/$\GL$ is incomplete w.r.t.\ $\tau_2$'' is consistent iff it
is consistent that Mahlo cardinals exist. \ec

It is possible to generalize these results to the case of bimodal
logic $\GLP_2$ \cite{Bek10b}. The situation remains essentially
unchanged, although a proof of Statement (i) of Blass's theorem
needs considerable adaptation.

\bt If $V=L$ and $\Omega\geq\aleph_\gw$, then $\GLP_2$ is complete
w.r.t.\ $(\Omega;\tau_1,\tau_2)$. \et

\section{Topologies for the variable-free fragment of $\GLP$}

A natural topological model for the variable-free fragment of $\GLP$
has been introduced by Thomas Icard \cite{Ica09}. It is not a
GLP-space, thus, not a model of the full $\GLP$. However, it is
sound and complete for the variable-free fragment of $\GLP$. It
gives a convenient tool for the study of this fragment, which plays
an important role in proof-theoretic applications of the polymodal
provability logic. Here, we give a simplified presentation of
Icard's polytopological space.

Let $\Omega$ be an ordinal and let $\ell:\Omega\to\Omega$ denote the
rank function for the order topology on $\Omega$. It is easy to
check that $\ell(0)=0$ and $\ell(\ga)=\gb$ whenever $\ga=\gy
+\gw^\gb$. By the Cantor normal form theorem, such a $\gb$ is
uniquely defined for each ordinal $\ga>0$. We define
$\ell^0(\ga)=\ga$ and $\ell^{k+1}(\ga)=\ell\ell^k(\ga)$.

Icard's topologies $\upsilon_n$, for each $n\in\gw$, are defined as
follows. Let $\upsilon_0$ be the left topology, and let $\upsilon_n$
be generated by $\upsilon_0$ and all sets of the form
$$U^m_\gb:=\{\ga\in \Omega:\ell^m(\ga)>\gb\},$$ for $m<n$.

Clearly, $\upsilon_n$ is an increasing sequence of topologies. In
fact, $\upsilon_1$ is the standard interval topology. We let $d_n$
and $\rho_n$ denote the derivative operator and the rank function
for $\upsilon_n$, respectively. We have the following
characterizations.

\bl \begin{enumr}
\item $\ell:(\Omega,\upsilon_{n+1})\to (\Omega,\upsilon_n)$ is a $d$-map;
\item $\upsilon_{n+1}$ is the coarsest topology $\nu$ on $\Omega$ such that
$\nu$ contains the interval topology and $\ell:(\Omega,\nu)\to
(\Omega,\upsilon_n)$ is continuous;
\item $\ell^n$ is the rank function of $\upsilon_n$, that is, $\rho_n=\ell^n$;
\item $\upsilon_{n+1}$ is generated by $\upsilon_n$ and $\{d_n^{\ga+1}(\Omega): \ga<\rho_n(\Omega)\}$.
\end{enumr}
\el

\bp (i) The map $\ell:(\Omega,\upsilon_{n+1})\to
(\Omega,\upsilon_n)$ is continuous. In fact, $\ell^{-1}[0,\gb)$ is
open in the order topology $\upsilon_1$, since
$\ell:(\Omega,\upsilon_1)\to (\Omega,\upsilon_0)$ is its rank
function, hence a $d$-map. Also, if $m<n$ then $\ell^{-1}(U^m_\gb)=
U^{m+1}_\gb$, hence it is open in $\upsilon_{n+1}$.

The map $\ell$ is open. Notice that $\upsilon_{n+1}$ is generated by
$\upsilon_1$ and some sets of the form $\ell^{-1}(U)$ where
$U\in\upsilon_n$. A base of $\upsilon_{n+1}$ consists of sets of the
form $V\cap \ell^{-1}(U)$, for some $V\in\upsilon_1$ and
$U\in\upsilon_n$. We have $\ell(V\cap \ell^{-1}(U))=\ell(V)\cap U$.
$\ell(V)$ is $\upsilon_0$-open, since $\ell:(\Omega,\upsilon_1)\to
(\Omega,\upsilon_0)$ is a $d$-map and $V\in\upsilon_1$. Hence, the
image of any basic open in $\upsilon_{n+1}$ is open in $\upsilon_n$.

The map $\ell$ is pointwise discrete, since $\ell^{-1}\{\ga\}$ is
discrete in the order topology $\upsilon_1$, hence in
$\upsilon_{n+1}$.

(ii) By (i), $\ell:(\Omega,\upsilon_{n+1})\to (\Omega,\upsilon_n)$
is continuous, hence $\nu\subseteq \upsilon_{n+1}$. On the other
hand, if $\ell:(\Omega,\nu)\to (\Omega,\upsilon_n)$ is continuous,
then $\ell^{-1}(U^m_\gb)\in\nu$, for each $m<n$. Therefore
$U^{m}_\gb\in\nu$, for all $m$ such that $1\leq m\leq n$. Since
$\nu$ also contains the interval topology, we have
$\upsilon_{n+1}\subseteq \nu$.

(iii) By (i), we have that $\rho_{n}\circ\ell$ is a $d$-map from
$(\Omega,\upsilon_{n+1})$ to $(\Omega,\upsilon_0)$. Hence, it
coincides with the rank function for $\upsilon_{n+1}$,
$\rho_{n+1}=\rho_n\circ \ell$. The claim follows by an easy
induction on $n$.

(iv) Notice that, by (iii), $$d_n^{\gb+1}(\Omega)=\{\ga\in\Omega:
\rho_n(\ga)>\gb\}=\{\ga\in\Omega:\ell^n(\ga)>\gb\}=U^n_\gb.$$
Obviously, $\upsilon_{n+1}$ is generated by $\upsilon_n$ and
$U^n_\gb$, for all $\gb$. Hence, the claim. \ep

We call \emph{Icard space} a polytopological space of the form
$(\Omega;\upsilon_0,\upsilon_1, \dots)$. Icard originally considered
just $\Omega=\ge_0$. We are going to give an alternative proof of
the following theorem  \cite{Ica09}.

\bt[Icard] \label{Icard} Let $\phi$ be a variable-free
$\GLP$-formula.
\begin{enumr}
\item If $\GLP\vdash\phi$, then $(\Omega;\upsilon_0,\upsilon_1, \dots)\models\phi$.
\item If $\Omega\geq\ge_0$ and $\GLP\nvdash\phi$, then $(\Omega;\upsilon_0,\upsilon_1, \dots)\nmodels\phi$.
\end{enumr}
\et

\bp Within this proof we abbreviate $(\Omega;\upsilon_0,\upsilon_1,
\dots)$ by $\Omega$. To prove part (i) we first remark that all
topologies $\upsilon_n$ are scattered, hence all axioms of $\GLP$
except for P1 are valid in $\Omega$. Moreover, $\Log(\Omega)$ is
closed under the inference rules of $\GLP$. Thus, we only have to
show that the variable-free instances of axiom P1 are valid in
$\Omega$. This is sufficient, because any derivation of a
variable-free formula in $\GLP$ can be replaced by a derivation in
which only the variable-free formulas occur (replace all the
variables by the constant $\top$).

Let $\phi$ be a variable-free formula. We denote by $\phi^*$ its
uniquely defined interpretation in $\Omega$. The validity of an
instance of P1 for $\phi$ amounts to the fact that $d_m(\phi^*)$ is
open in $\upsilon_{n}$ whenever $m<n$. Thus, we have to prove the
following proposition.

\bpr \label{open} For any variable-free formula $\phi$,
$d_n(\phi^*)$ is open in $\upsilon_{n+1}$. \epr

Let $\phi^+$ denote the result of replacing in $\phi$ each modality
$\la n\ra$ by $\la n+1\ra$. We need the following auxiliary claim.

\bl \label{plus} If $\phi$ is variable-free, then
$\ell^{-1}(\phi^*)=(\phi^+)^*$. \el

\bp This goes by induction on the build-up of $\phi$. The cases of
constants and boolean connectives are easy. Suppose $\phi=\la n\ra
\psi$. We notice that since $\ell:(\Omega,\upsilon_{n+1})\to
(\Omega,\upsilon_n)$ is a $d$-map, we have
$\ell^{-1}(d_n(A))=d_{n+1}(\ell^{-1}(A))$, for any
$A\subseteq\Omega$. Therefore,
$$\ell^{-1}(\phi^*)=\ell^{-1}(d_n(\psi^*))=
d_{n+1}(\ell^{-1}(\psi^*))=d_{n+1}((\psi^+)^*)= (\phi^+)^*,$$ as
required. \ep

We prove Proposition \ref{open} in two steps. First, we show that it
holds for a subclass of variable-free formulas called \emph{ordered
formulas}. Then we show that any variable-free formula is equivalent
in $\Omega$ to an ordered one.

A formula $\phi$ is called \emph{ordered}, if no modality $\la m\ra$
occurs within the scope of $\la n\ra$ in $\phi$ for any $m<n$. The
\emph{height of $\phi$} is the index of its maximal modality.

\bl \label{order} If $\la n\ra\phi$ is ordered, then $d_n(\phi^*)$
is open in $\upsilon_{n+1}$. \el

\bp This goes by induction on the height of $\la n\ra\phi$. If it is
$0$, then $n=0$. If $n=0$, the claim is obvious, since $d_0(A)$ is
open in $\upsilon_1$, for any $A\subseteq\Omega$. If $n>0$, since
$\la n\ra\phi$ is ordered, we observe that $\la n\ra\phi$ has the
form $(\la n-1\ra\psi)^+$, for some $\psi$. The height of $\la
n-1\ra\psi$ is less than that of $\la n\ra\phi$. Hence, by the
induction hypothesis, $(\la n-1\ra\psi)^*\in\upsilon_n$. Since
$\ell:(\Omega,\upsilon_{n+1})\to (\Omega,\upsilon_n)$ is continuous,
we conclude that $\ell^{-1}(\la n-1\ra\psi)^*$ is open in
$\upsilon_{n+1}$. By Lemma \ref{plus}, this set coincides with $(\la
n\ra\phi)^*=d_n(\phi^*)$. \ep

\bl Any variable-free formula $\phi$ is equivalent in $\Omega$ to an
ordered one. \el \bp We argue by induction on the complexity of
$\phi$. The cases of boolean connectives and constants are easy.
Suppose $\phi$ has the form $\la n\ra\psi$, where we can assume
$\psi$ to be in disjunctive normal form
$\psi=\bigvee_i\bigwedge_j\pm\la n_{ij}\ra\psi_{ij}$. By the
induction hypothesis, we may assume all the subformulas $\la
n_{ij}\ra\psi_{ij}$ (and $\psi$ itself) ordered. Since $\la n\ra$
commutes with disjunction, it will be sufficient to show that, for
each $i$, the formula
$$\la n\ra \bigwedge_j\pm\la n_{ij}\ra\psi_{ij} \eqno (*)$$ can be ordered.

By Lemma \ref{order} each set $(\la n_{ij}\ra\psi_{ij})^*$ is open
in $\upsilon_n$, whenever $n_{ij}<n$. Being a derived set, it is
also closed in $\upsilon_{n_{ij}}$ and hence in $\upsilon_n$. Thus,
all such sets are clopen.

If $U$ is open then $d(A\cap U)=d (A) \cap U$, for any topological
space. In particular, for any $A\subseteq \Omega$ and $n_{ij}<n$,
$$d_n(A\cap (\pm\la n_{ij}\ra\psi_{ij})^*)= d_n (A) \cap (\pm\la n_{ij}\ra\psi_{ij})^*.$$
This allows us to bring all the conjuncts $\pm\la
n_{ij}\ra\psi_{ij}$ from under the $\la n\ra$ modality in $(*)$. The
resulting conjunction is ordered. \ep

This concludes the proof of Proposition \ref{open} and thereby of
Part (i).

\bigskip
To prove Part (ii), we shall rely on the following fundamental lemma
about the variable-free fragment of $\GLP$. For a purely syntactic
proof of this lemma we refer the reader to \cite{Bek04,Bek05a}.

A variable-free formula $A$ is called a \emph{word} if it is
built-up from $\top$ only using connectives of the form $\la n\ra$,
for any $n\in\gw$. We write $A\vdash B$ for $\GLP\vdash A\to B$.

\bl \label{words}
\begin{enumr}
\item Every variable-free formula is equivalent in $\GLP$ to a boolean combination of words;
\item For any words $A$ and $B$, either $A\vdash \la 0 \ra B$, or $B\vdash \la 0 \ra A$, or $A$ and $B$ are equivalent;
\item Conjunction of words is equivalent to a word.
\end{enumr}
\el

We prove Part (ii) of Theorem \ref{Icard} in a series of lemmas.
First, we show that any word is satisfiable in $\Omega$ provided
$\Omega\geq \ge_0$.

\bl For any word $A$,  $\ge_0\in A^*$. \el \bp We know that
$\rho_n(\ge_0)=\ell^n(\ge_0)=\ge_0$. Hence, $\ge_0\in d_n(\Omega)$,
for each $n$. Assume $n$ exceeds all the indices of modalities in
$A$ and $A=\la m\ra B$. By Proposition \ref{open} the set $B^*$ is
open in $\upsilon_n$. By the induction hypothesis $\ge_0\in B^*$.
Hence, $\ge_0\in d_n(B^*)\subseteq d_m(B^*)= A^*$. This proves the
claim. \ep

Applying this lemma to the word $\la 0\ra A$ we obtain the following
corollary.

\bc For every word $A$, there is an $\ga<\ge_0$ such that $\ga\in
A^*$. \ec

Let $\min(A^*)$ denote the least ordinal $\ga\in \Omega$ such that
$\ga\in A^*$.

\bl \label{min} For any words $A,B$, if  $A\nvdash B$ then
$\min(A^*)\notin B^*$. \el

\bp If $A\nvdash B$ then, by Lemma \ref{words} (ii), $B\vdash \la
0\ra A$. Therefore, by the soundness of $\GLP$ in $\Omega$,
$B^*\subseteq d_0(A^*)$. It follows that for each $\gb\in B^*$ there
is an $\ga\in A^*$ such that $\ga<\gb$. Thus, $\min(A^*)\notin B^*$.
\ep

Now we are ready to prove Part (ii). Assume $\phi$ is variable-free
and $\GLP\nvdash \phi$. By Lemma \ref{words} (i) we can assume that
$\phi$ is a boolean combination of words. Writing $\phi$ in
conjunctive normal form we observe that it is sufficient to prove
the claim only for formulas $\phi$ of the form $\bigwedge_i A_i \to
\bigvee_j B_j,$ where $A_i$ and $B_j$ are words. Moreover,
$\bigwedge_i A_i$ is equivalent to a single word $A$.

Since $\GLP\nvdash \phi$ we have $A\nvdash B_j$, for each $j$. Let
$\ga=\min(A^*)$. By Lemma \ref{min} we have $\ga \notin B_j^*$, for
each $j$. Hence, $\ga\notin (\bigvee_j B_j)^*$ and $\ga\notin
\phi^*$. This means that $\Omega\nmodels\phi^*$. \ep

\section{Further results}

Topological semantics of polymodal provability logic has been
extended to the language with transfinitely many modalities. A
logic $\GLP_\Lambda$ having modalities $[\ga]$, for all ordinals
$\ga<\Lambda$, is introduced in \cite{Bek05a}. It was intended for
the proof-theoretic analysis of predicative theories and is currently being actively investigated for that purpose.

David Fernandez and Joost Joosten undertook a thorough study of the
variable-free fragment of that logic mostly in connection with the
arising ordinal notation systems (see~\cite{JF12-aiml,JF12-Turing}
for a sample). In particular, they found a suitable generalization
of Icard's polytopological space and showed that it is complete for
that fragment \cite{JF12-models}. Fernandez \cite{Fern12-top} also
proved topological completeness of the full $\GLP_\Lambda$ by
generalizing the results in \cite{BekGab11}.

The ordinal GLP-space is easily generalized to transfinitely many
topologies $(\tau_\ga)_{\ga<\Lambda}$ by letting $\tau_0$ be the
left topology, $\tau_{\ga+1}:=\tau_\ga^+$ and, for limit ordinals
$\gl$, $\tau_\gl$ be the topology generated by all $\tau_\ga$ such
that $\ga<\gl$. This space is a natural model of $\GLP_\Lambda$ and
has been studied quite recently by Joan Bagaria \cite{Bag12} and
further by Bagaria, Magidor and Sakai \cite{BaMaSa12}. In
particular, the three authors proved that in $L$ the limit points of
$\tau_{n+2}$ are $\Pi_n^1$-indescribable cardinals. The question
posed in \cite{BBI09} whether the non-discreteness of $\tau_{n+2}$
is equiconsistent with the existence of $\Pi_n^1$-indescribable
cardinals still seems to be open.

\section{Acknowledgement} The first author was supported by the
Russian Foundation for Basic Research (RFBR), Russian Presidential
Council for Support of Leading Scientific Schools, and the
Swiss--Russian cooperation project STCP--CH--RU ``Computational
proof theory''.
\newline\indent The second author was supported by the Shota
Rustaveli National Science Foundation grant {\#}FR/489/5-105/11 and
the French--Georgian grant CNRS--SRNSF {\#}4135/05-01.


\begin{thebibliography}{10}

\bibitem{Aba85}
M.~Abashidze.
\newblock Ordinal completeness of the {G}\"{o}del-{L}\"{o}b modal system.
\newblock In {\em Intensional logics and the logical structure of theories},
  pages 49--73. Metsniereba, Tbilisi, 1985.
\newblock In Russian.

\bibitem{AB04}
S.N. Artemov and L.D. Beklemishev.
\newblock Provability logic.
\newblock In D.~Gabbay and F.~Guenthner, editors, {\em {Handbook of
  Philosophical Logic}, 2nd ed.}, volume~13, pages 229--403. Kluwer, Dordrecht,
  2004.

\bibitem{Bag12}
J.~Bagaria.
\newblock Topologies on ordinals and the completeness of polymodal provability
  logic.
\newblock In preparation, 2012.

\bibitem{BaMaSa12}
J.~Bagaria, M.~Magidor, and H.~Sakai.
\newblock Reflection and indescribability in the constructible universe.
\newblock Manuscript, 2012.

\bibitem{BekGab11}
L.~Beklemishev and D.~Gabelaia.
\newblock Topological completeness of the provability logic {GLP}.
\newblock Preprint arXiv:1106.5693v1 [math.LO], 2011.

\bibitem{Bek04}
L.D. Beklemishev.
\newblock Provability algebras and proof-theoretic ordinals, {I}.
\newblock {\em {Annals of Pure and Applied Logic}}, 128:103--123, 2004.

\bibitem{Bek05}
L.D. Beklemishev.
\newblock Reflection principles and provability algebras in formal arithmetic.
\newblock {\em {Uspekhi Matematicheskikh Nauk}}, 60(2):3--78, 2005.
\newblock In Russian. English translation in: \emph{Russian Mathematical
  Surveys}, 60(2): 197--268, 2005.

\bibitem{Bek05a}
L.D. Beklemishev.
\newblock Veblen hierarchy in the context of provability algebras.
\newblock In P.~H\'ajek, L.~Vald\'es-Villanueva, and D.~Westerst{\aa}hl,
  editors, {\em {Logic, Methodology and Philosophy of Science, Proceedings of
  the Twelfth International Congress}}, pages 65--78. Kings College
  Publications, London, 2005.
\newblock Preprint: Logic Group Preprint Series 232, Utrecht University, June
  2004.

\bibitem{Bek06}
L.D. Beklemishev.
\newblock The {Worm Principle}.
\newblock In Z.~Chatzidakis, P.~Koepke, and W.~Pohlers, editors, {\em {Lecture
  Notes in Logic 27. Logic Colloquium '02}}, pages 75--95. AK Peters, 2006.
\newblock Preprint: Logic Group Preprint Series 219, Utrecht University, March
  2003.

\bibitem{Bek09a}
L.D. Beklemishev.
\newblock On {GLP}-spaces.
\newblock Manuscript, http://www.mi.ras.ru/~bekl/Papers/glp-sp.pdf, 2009.

\bibitem{Bek10}
L.D. Beklemishev.
\newblock Kripke semantics for provability logic {GLP}.
\newblock {\em {Annals of Pure and Applied Logic}}, 161:756--774, 2010.
\newblock Preprint: Logic Group Preprint Series 260, University of Utrecht,
  November 2007. \texttt{http://preprints.phil.uu.nl/lgps/}.

\bibitem{Bek10b}
L.D. Beklemishev.
\newblock Ordinal completeness of bimodal provability logic {GLB}.
\newblock In N.~Bezhanishvili et~al., editor, {\em Logic, Language and
  Computation, TbiLLC 2009. Lecture Notes in Artificial Intelligence}, number
  6618, pages 1--15, 2011.

\bibitem{Bek11}
L.D. Beklemishev.
\newblock A simplified proof of the arithmetical completeness theorem for the
  provability logic {GLP}.
\newblock {\em Trudy Matematicheskogo Instituta imeni V.A.Steklova},
  274(3):32--40, 2011.
\newblock In Russian. English translation: \emph{Proceedings of the Steklov
  Institute of Mathematics}, 274(3):25--33, 2011.

\bibitem{BBI09}
L.D. Beklemishev, G.~Bezhanishvili, and T.~Icard.
\newblock On topological models of {GLP}.
\newblock In R.~Schindler, editor, {\em Ways of Proof Theory}, Ontos
  Mathematical Logic, pages 133--152. Ontos Verlag, Heusendamm bei Frankfurt,
  Germany, 2010.
\newblock Preprint: Logic Group Preprint Series 278, University of Utrecht,
  August 2009, \texttt{http://preprints.phil.uu.nl/lgps/}.

\bibitem{BJV}
L.D. Beklemishev, J.~Joosten, and M.~Vervoort.
\newblock A finitary treatment of the closed fragment of {J}aparidze's
  provability logic.
\newblock {\em Journal of Logic and Computation}, 15(4):447--463, 2005.

\bibitem{BEG05}
G.~Bezhanishvili, L.~Esakia, and D.~Gabelaia.
\newblock Some results on modal axiomatization and definability for topological
  spaces.
\newblock {\em Studia Logica}, 81:325--355, 2005.

\bibitem{BM10}
G.~Bezhanishvili and P.~Morandi.
\newblock Scattered and hereditarily irresolvable spaces in modal logic.
\newblock {\em Archive}, 49(3):343--365, 2010.

\bibitem{Bla90}
A.~Blass.
\newblock Infinitary combinatorics and modal logic.
\newblock {\em Journal of Symbolic Logic}, 55(2):761--778, 1990.

\bibitem{Boo79}
G.~Boolos.
\newblock {\em The Unprovability of Consistency: An Essay in Modal Logic}.
\newblock Cambridge University Press, Cambridge, 1979.

\bibitem{Boo93a}
G.~Boolos.
\newblock The analytical completeness of {D}zhaparidze's polymodal logics.
\newblock {\em Annals of Pure and Applied Logic}, 61:95--111, 1993.

\bibitem{Boo93}
G.~Boolos.
\newblock {\em The {L}ogic of {P}rovability}.
\newblock Cambridge University Press, Cambridge, 1993.

\bibitem{Car05}
L.~Carlucci.
\newblock Worms, gaps and hydras.
\newblock {\em Mathematical Logic Quarterly}, 51(4):342--350, 2005.

\bibitem{ChZa}
A.~Chagrov and M.~Zakharyaschev.
\newblock {\em Modal Logic}.
\newblock Clarendon Press, Oxford, 1997.

\bibitem{JF12-aiml}
{D. Fern\'andez-Duque and J. J. Joosten}.
\newblock Kripke models of transfinite provability logic.
\newblock In {\em Advances in Modal Logic}, volume~9, pages 185--199. College
  Publications, 2012.

\bibitem{JF12-models}
{D. Fern\'andez-Duque and J. J. Joosten}.
\newblock Models of transfinite provability logic.
\newblock To appear in \emph{Journal of Symbolic Logic}. Preliminary version at
  http://personal.us.es/dfduque/glpmodels.pdf, 2012.

\bibitem{JF12-Turing}
{D. Fern\'andez-Duque and J. J. Joosten}.
\newblock Turing progressions and their well-orders.
\newblock In {\em How the world computes}, Lecture Notes in Computer Science,
  pages 212--221. Springer, 2012.

\bibitem{Esa81}
L.~Esakia.
\newblock Diagonal constructions, {L}\"{o}b's formula and {C}antor's scattered
  spaces.
\newblock In {\em Studies in logic and semantics}, pages 128--143. Metsniereba,
  Tbilisi, 1981.
\newblock In Russian.

\bibitem{God-CW-I}
S.~Feferman, J.R. Dawson, S.C. Kleene, G.H. Moore, R.M. Solovay, and J.~van
  Heijenoort, editors.
\newblock {\em Kurt G\"{o}del Collected Works, Volume 1: Publications
  1929--1936}. Oxford Univeristy Press, 1996.

\bibitem{Fern12-top}
D.~Fern\'andez-Duque.
\newblock The polytopologies of transfinite provability logic.
\newblock Preprint arXiv:1207.6595v2 [math.LO], 2012.

\bibitem{God33}
K.~G\"{o}del.
\newblock Eine {I}nterpretation des intuitionistischen {A}ussagenkalkuls.
\newblock {\em Ergebnisse Math. Kolloq.}, 4:39--40, 1933.
\newblock English translation in \cite{God-CW-I}, pages 301--303.

\bibitem{HarShe}
L.~Harrington and S.~Shelah.
\newblock Some exact equiconsistency results in set theory.
\newblock {\em Notre Dame Journal of Formal Logic}, 26:178--188, 1985.

\bibitem{IcaJoo12}
T.~Icard and J.~Joosten.
\newblock Provability and interpretability logics with restricted realizations.
\newblock {\em Notre Dame Journal of Formal Logic}, 53(2):133--154, 2012.

\bibitem{Ica09}
T.F. Icard, III.
\newblock A topological study of the closed fragment of {GLP}.
\newblock {\em Journal of Logic and Computation}, 21(4):683--696, 2011.

\bibitem{Ign93}
K.N. Ignatiev.
\newblock On strong provability predicates and the associated modal logics.
\newblock {\em The Journal of Symbolic Logic}, 58:249--290, 1993.

\bibitem{Dzh86}
G.K. Japaridze.
\newblock The modal logical means of investigation of provability.
\newblock {T}hesis in {P}hilosophy, in {R}ussian, Moscow, 1986.

\bibitem{Jech}
T.~Jech.
\newblock {\em Set Theory. The Third Millenium Edition}.
\newblock Springer, 2002.

\bibitem{Kan09}
A.~Kanamori.
\newblock {\em The Higher Infinite. Second edition}.
\newblock Springer, 2009.

\bibitem{Lob55}
M.H. L\"{o}b.
\newblock Solution of a problem of {Leon Henkin}.
\newblock {\em The Journal of Symbolic Logic}, 20:115--118, 1955.

\bibitem{MS73}
A.~Macintyre and H.~Simmons.
\newblock G\"odel's diagonalization technique and related properties of
  theories.
\newblock {\em Colloq. Math.}, 28:165--180, 1973.

\bibitem{Mag}
R.~Magari.
\newblock The diagonalizable algebras (the algebraization of the theories which
  express {T}heor.:{II}).
\newblock {\em Bollettino della Unione Matematica Italiana, {\rm Serie 4}}, 12,
  1975.
\newblock Suppl. fasc. 3, 117--125.

\bibitem{Magi82}
M.~Magidor.
\newblock Reflecting stationary sets.
\newblock {\em The Journal of Symbolic Logic}, 47:755--771, 1982.

\bibitem{MT44}
J.~C.~C. McKinsey and A.~Tarski.
\newblock The algebra of topology.
\newblock {\em Annals of Mathematics}, 45:141--191, 1944.

\bibitem{MekShe}
A.~Mekler and S.~Shehlah.
\newblock The consistency strength of ``every stationary set reflects''.
\newblock {\em Israel Journal of Mathematics}, 67(3):353--366, 1989.

\bibitem{Sim75}
H.~Simmons.
\newblock Topological aspects of suitable theories.
\newblock {\em Proc. Edinburgh Math. Soc.~(2)}, 19:383--391, 1975.

\bibitem{Smo85}
C.~Smory\'{n}ski.
\newblock {\em Self-Reference and Modal Logic}.
\newblock Springer-Verlag, Berlin, 1985.

\bibitem{Sol76}
R.M. Solovay.
\newblock Provability interpretations of modal logic.
\newblock {\em Israel Journal of Mathematics}, 28:33--71, 1976.

\bibitem{Ste07}
C.~Steinsvold.
\newblock {\em Topological models of belief logics}.
\newblock PhD thesis, New York, 2007.

\bibitem{Vis84}
A.~Visser.
\newblock The provability logics of recursively enumerable theories extending
  {P}eano {A}rithmetic at arbitrary theories extending {P}eano {A}rithmetic.
\newblock {\em Journal of Philosophic Logic}, 13:97--113, 1984.

\end{thebibliography}

\end{document}